\documentclass[review]{elsarticle}
\usepackage{pifont}
\usepackage[top=1in, left=1in, bottom=1in, right=1in]{geometry}
\usepackage{fleqn}
\usepackage{hyperref}
\usepackage{graphicx}
\usepackage{amsmath}
\usepackage{amsfonts}
\usepackage{amssymb}
\usepackage{natbib}
\usepackage{color}
\usepackage{float}
\newtheorem{theorem}{Theorem}

\newtheorem{property}[theorem]{Property}

\newenvironment{proof}[1][Proof]{\noindent\textbf{#1.} }{\ \rule{0.5em}{0.5em}}
\newenvironment{remark}[1][Remark]{\noindent\textbf{#1.} }%{\ \rule{0.5em}{0.5em}}

\begin{document}

\title{The Role of Constraints in a Segregation Model: The Symmetric Case}

\author{D. Radi\fnref{fn1}}
\ead{d.radi@univpm.it}
\address[fn1]{Department of Management, Polytechnic University of Marche, Ancona, Italy}

\author{L. Gardini\corref{cor1}\fnref{fn2}}
\ead{laura.gardini@uniurb.it}
\cortext[cor1]{Corresponding author}
\address[fn2]{DESP, University of Urbino "Carlo Bo", Urbino, Italy}

\author{V. Avrutin\fnref{fn3,fn2}}
\ead{Viktor.Avrutin@ist.uni-stuttgart.de}
\address[fn3]{IST, University of Stuttgart, Stuttgart, Germany}

\begin{abstract}

In this paper we study the effects of constraints on the dynamics of an
adaptive segregation model introduced by Bischi and Merlone (2011). The model
is described by a two dimensional piecewise smooth dynamical system in
discrete time. It models the dynamics of entry and exit of two populations
into a system, whose members have a limited tolerance about the presence of
individuals of the other group. The constraints are given by the upper limits
for the number of individuals of a population that are allowed to enter the
system. They represent possible exogenous controls imposed by an authority in
order to regulate the system. Using analytical, geometric and numerical
methods, we investigate the border collision bifurcations generated by these
constraints assuming that the two groups have similar characteristics and have
the same level of tolerance toward the members of the other group. We also
discuss the policy implications of the constraints to avoid
segregation.\bigskip
\end{abstract}

%\textbf{Keywords:}
\maketitle

\section{Introduction}\label{Intro}

In his seminal contribution \cite{Schelling1969}, Schelling underlines how
discriminatory individual choices can lead to the segregation of two groups of
people of opposite kind. People get separated for different reasons, such as
sex, age, income, language or nationality, color of the skin, and the like.
Since then, this idea has been developed and tested using mainly an agent
based computer simulation approach, see e.g., \cite{EpsteinAxtell1996} and
\cite{Zhang2004}. Instead, \cite{BischiMerlone2011} introduces an adaptive
dynamical model in discrete time that captures the features of the segregation
process designed by Schelling. This model is represented by an iterated two
dimensional non invertible map. The analysis of the model provides a rather
solid mathematical ground that confirms and extends the qualitative
illustration of the dynamics provided by Schelling in \cite{Schelling1969}. In
particular, the possibility, depending on the initial conditions, to end up
either in an equilibrium of segregation or an equilibrium of coexistence of
the members of the two groups in the same system. The investigation reveals
also more complicated phenomena which could have not been observed in
\cite{Schelling1969} due to the lack of mathematical formalization of the model,
such as the emergence of periodic or chaotic solutions. Such oscillatory
solutions represent situations in which the number of the members of the two
groups that enter or exit the system oscillate perpetually in time as the
results of overshooting due to impulsive (or emotional) behavior of the agents.

Following Schelling's ideas, the authors of \cite{BischiMerlone2011}
introduced in the model two constraints that limit the maximum number of the
members of each group allowed to enter the system. This is indeed quite
relevant, as the constraint may reflect the policy decision of some state or
group. The constraints make the model piecewise differentiable and, from a
dynamical point of view, can be responsible for possible border collision
bifurcations. In \cite{BischiMerlone2011}, the effects of these constraints
are only marginally analyzed and a deeper investigation is left for further
researches. In this paper, following their suggestion we provide a
comprehensive description of the effects of these constraints on the dynamics
of the model. In particular, we use geometrical, analytical and numerical
tools to investigate the nature of the dynamics that can arise changing the
value of these constraints.

Limiting the analysis to a symmetric setting, i.e. assuming that the two
populations are of the same size and have the same level of tolerance toward
the other type of agents, it emerges that if the two constraints are both
sufficiently tight, then an equilibrium of non segregation exists and it is
stable, together with two coexisting equilibria of segregation, which are
always present and always stable. In particular, the two-dimensional
bifurcation diagram reveals that if we relax the limitations to the maximum
number of the members of the two populations allowed to enter the system, then for certain initial conditions
we first observe a transition from a stable equilibrium of coexistence to
stable cycles of any periodicity and subsequently a transition from stable
cycles to equilibria of segregation. On the contrary, if the
constraints are not fixed equally, for example we limit more the members of
the population one to enter the system and less the members of population two
and this gap is large enough, as a result we can have either only stable
equilibria of segregation or coexistence of a stable periodic solution and
stable equilibria of segregation. Thus, it is necessary to impose equal and sufficiently tight constraints on the maximum number of the members of the two
populations allowed to enter the system, to have, at least for certain initial conditions, the possibility to convergence to an equilibrium of non segregation.

The dynamics of the model here proposed are particularly interesting from a
mathematical point of view as well. Indeed, the model is described by a
continuous two-dimensional piecewise differentiable map, with several borders
crossing which the system changes its definition. The dynamics associated with
piecewise smooth systems is a quite new research branch, and several papers
have been dedicated to this subject in the last decade (see, e.g., \cite{Bernardo08} and \cite{ZhusubaliyevMosekilde2003}). Such an increasing interest towards nonsmooth dynamics comes both
from the new theoretical problems due to the borders and from the wide
interest in the applied context. In fact, many models are described by
constrained functions, leading to piecewise smooth systems, continuous or
discontinuous. We recall several oligopoly models with different kinds of
constraints considered in the books \cite{PuuSushko2002} and \cite{BischiChiarellaKopelSzidarovszky2009}, nonsmooth
economic models in \cite{Day1994}, \cite{HommesNusse1991}, \cite{Matsuyama2010}, \cite{PuuSushko2006} and
\cite{GardiniSushkoNaimzada2008}, financial market modeling in \cite{DayChen1993}, \cite{TramontanaWesterhoffGardini2010} and
\cite{TramontanaGardiniWesterhoff2011}, and modeling of multiple-choice in \cite{BischiGardiniMerlone2009}, \cite{GardiniMerloneTramontana2011} and
\cite{DalFornoGardiniMerlone2012}.

The map considered in the present paper is characterized by several
constraints, leading to several different partitions of the phase plane in
which the system changes definition. Moreover, the definitions in some regions
are quite degenerate, as mapped into points or segments of straight lines.
That is, the degeneracy consists in a Jacobian matrix which has one or two
eigenvalues equal to zero in the points of a whole region. Thus, when an
invariant set as a cycle has a periodic point colliding with a border, then a
border collision occurs, which often leads to a \textit{border collision
bifurcation} (BCB for short), first described in \cite{NusseYorke1992} (see also
\cite{NusseYorke1995} and \cite{SushkoGardini2010}). The result of the contact, that is, what happens
to the dynamics after the contact, is in general difficult to predict.
However, in one-dimensional piecewise smooth systems, the possible results of
a generic BCB of an attracting cycle with one border point can be rigorously
classified depending on the parameters using the one-dimensional\textit{\ BCB
normal form,} which is the well known skew tent map defined by two linear
functions. In fact, the dynamics of the skew tent map are completely described
according to the slopes of the linear branches, and it is possible to use this
map as a normal form (see, e.g., \cite{ItoTanakaNakada1979}, \cite{MastrenkoMastrenkoChua1993}, \cite{SushkoAgliariGardini2006} and
\cite{SushkoGardini2010}).

This powerful result will be used also in the analysis of the two-dimensional
system considered in this work. This is due to the high degeneracy of the map,
often leading to a dynamic behavior which is constrained to some
one-dimensional set, and in it the map can be studied by using its
one-dimensional restriction. Another peculiarity of the degeneracy (when the
system is defined by constant values in one or both variables), is that the
one-dimensional restriction is characterized by a flat branch in the shape of
the function. For a piecewise smooth map with a flat branch any cycle with a
point on that branch is \textit{superstable} (i.e. it has a $0$ eigenvalue).
Moreover, in the applied context it is important to stress that superstable
cycles related to a flat branch, differently from "smooth" superstable cycles,
are persistent under parameters' perturbations. That is, in the parameter
space there are open regions related to these cycles, as we shall see also in
our map. Clearly, the boundaries of such periodicity regions can be defined
only by BCBs of the related cycles given that the zero eigenvalue doesn't
allow any other bifurcation. Examples of systems characterized by a map with a
flat branch can be found in \cite{BrianzoniMichettiSushko2010}, \cite{AvrutinFutterSchanz2012}, \cite{TramontanaGardiniPuu2011} and
\cite{SushkoGardiniMatsuyama2014}. The feature of such systems is that the bifurcation structure in
some of the periodicity regions of superstable cycles of the parameter space
are organized according to the well known \textit{U-sequence} (first described
in \cite{MetropolisSteinStein1973}, see also \cite{Hao1989}) which is characteristic for unimodal maps.
In \cite{SushkoGardiniMatsuyama2014} this is well described introducing one more letter related to
the flat branch, besides the two-letters for the symbolic sequences in
increasing/decreasing branches. In the U-sequence the BCB are related to
infinite cascades of\textit{\ flip BCBs }(not standard flip, as not related to
eigenvalues), and the first symbolic sequence in such a cascade for the cycles
of periods $n>1$ is related to the cycle born due to \textit{fold BCB }(not
standard fold, or tangent, bifurcation as in smooth maps).

The plan of the work is as follows. In Section \ref{ModSetup} we introduce the model and
describe its main dynamical properties. In Section \ref{BCBs} we analyze the effect of
the constraints on the dynamics of the model. In particular, we investigate
the BCBs that occur as the constraints change and we provide the main
implications in terms of segregation. In Section \ref{Conc}, we conclude providing some
indications for possible further explorations of the dynamics of the
model.\bigskip

\section{Model setup and preliminaries}\label{ModSetup}

As in \cite{BischiMerlone2011} and \cite{Schelling1969}, we assume that
individuals are partitioned in two classes $C_{1}$ and $C_{2}$, say "group 1"
and "group 2", of respective numerosity $N_{1}$ and $N_{2}$ and that each
group cares about the type of the people in the district they live in.

Moreover, we assume that any individual of group $i,\;i=1,2$, can observe the
ratio of the two types of agents at any moment, and can decide to move in
(out) depending on its own (dis)satisfaction with the observed proportion of
opposite type agents to its own type. This degree of (dis)satisfaction is
capture by functions%
\begin{equation}
R_{1}\left(  x_{1}\right)  =\tau_{1}\left(  1-\frac{x_{1}}{N_{1}}\right)
\ \text{and }R_{2}\left(  x_{2}\right)  =\tau_{2}\left(  1-\frac{x_{2}}{N_{2}%
}\right)
\end{equation}
where $x_{i}R\left(  x_{i}\right)  $ gives the maximum number of agents of
group $C_{j}$, that are tolerated by $x_{i}$ agents of group $C_{i}$. It
follows that agents of type $i$ will enter the system if $x_{i}R_{i}\left(
x_{i}\right)  -x_{j}>0$ and will exit otherwise. From which we have that the
equation giving the number of agents of type $i$ that are in the system at
time $t+1$ is%
\begin{equation}
\frac{x_{i}\left(  t+1\right)  -x_{i}\left(  t\right)  }{x_{i}\left(
t\right)  }=\gamma_{i}\left[  x_{i}\left(  t\right)  R_{i}\left(  x_{i}\left(
t\right)  \right)  -x_{j}\left(  t\right)  \right]
\end{equation}
where $\gamma_{i}$ is the speed of adjustment. Assuming also a restriction on
the number of members of group $C_{i}$ that are allowed to enter the system,
say $0\leq x_{i}\left(  t\right)  \leq K_{i}$, with $K_{i}\leq N_{i}$, as a
result we obtain the following segregation model, as proposed in
\cite{BischiMerlone2011}, which is rich of different dynamic behaviors. It is
described by a continuous two-dimensional piecewise-smooth map $T:\ R_{+}%
^{2}\rightarrow R_{+}^{2}$ given by%

\begin{equation}
(x_{1}(t+1),x_{2}(t+1))=T(x_{1}(t),x_{2}(t))=(T_{1}(x_{1}(t),x_{2}%
(t)),T_{2}(x_{1}(t),x_{2}(t)))\label{T}%
\end{equation}
with%
\begin{equation}
T_{1}\left(  x_{1},x_{2}\right)  =\left\{
\begin{array}
[c]{ccc}%
0 & \text{if} & F_{1}\left(  x_{1},x_{2}\right)  \leq0\\
F_{1}\left(  x_{1},x_{2}\right)  & \text{if} & 0\leq F_{1}\left(  x_{1}%
,x_{2}\right)  \leq K_{1}\\
K_{1} & \text{if} & F_{1}\left(  x_{1},x_{2}\right)  \geq K_{1}%
\end{array}
\right.
\end{equation}%
\begin{equation}
T_{2}\left(  x_{1},x_{2}\right)  =\left\{
\begin{array}
[c]{ccc}%
0 & \text{if} & F_{2}\left(  x_{1},x_{2}\right)  \leq0\\
F_{2}\left(  x_{1},x_{2}\right)  & \text{if} & 0\leq F_{2}\left(  x_{1}%
,x_{2}\right)  \leq K_{2}\\
K_{2} & \text{if} & F_{2}\left(  x_{1},x_{2}\right)  \geq K_{2}%
\end{array}
\right.
\end{equation}
where%
\begin{align}
F_{1}\left(  x_{1},x_{2}\right)   & =x_{1}\left[  1-\gamma_{1}x_{2}+\gamma
_{1}x_{1}R_{1}\left(  x_{1}\right)  \right] \\
F_{2}\left(  x_{1},x_{2}\right)   & =x_{2}\left[  1-\gamma_{2}x_{1}+\gamma
_{2}x_{2}R_{2}\left(  x_{2}\right)  \right] \nonumber
\end{align}
Let us also recall the conditions on the parameters. We have that for $i=1,2,
$ $\gamma_{i},$ $\tau_{i}$ and $N_{i}$ can take any positive value, and it
must be $K_{i}\leq N_{i}.$

From the definition of the map we have that the phase plane of the dynamical
system can be divided into several regions where the system is defined by
different functions. On the boundaries of the regions the map is continuous
but not differentiable. The boundaries of non differentiability are given by
the curves $F_{i}\left(  x_{1},x_{2}\right)  =K_{i}$\ which can be written in
explicit form as follows:%

\begin{equation}%
\begin{array}
[c]{lll}%
BC_{1,K}:x_{2}=\left[  1+\gamma_{1}x_{1}R_{1}\left(  x_{1}\right)
-\frac{K_{1}}{x_{1}}\right]  /\gamma_{1} & \text{where} & F_{1}\left(
x_{1},x_{2}\right)  =K_{1}\\
BC_{2,K}:x_{1}=\left[  1+\gamma_{2}x_{2}R_{2}\left(  x_{2}\right)
-\frac{K_{2}}{x_{2}}\right]  /\gamma_{2} & \text{where} & F_{2}\left(
x_{1},x_{2}\right)  =K_{2}%
\end{array}
\label{BCK}%
\end{equation}
and $F_{i}\left(  x_{1},x_{2}\right)  =0$ which, as it is immediate, are
satisfied by $x_{i}=0,$ and other points belonging to the curves given by:%
\begin{equation}%
\begin{array}
[c]{lll}%
BC_{1,0}:x_{2}=\left[  1+\gamma_{1}x_{1}R_{1}\left(  x_{1}\right)  \right]
/\gamma_{1} & \text{where} & F_{1}\left(  x_{1},x_{2}\right)  =0,~x_{1}\neq0\\
BC_{2,0}:x_{1}=\left[  1+\gamma_{2}x_{2}R_{2}\left(  x_{2}\right)  \right]
/\gamma_{2} & \text{where} & F_{2}\left(  x_{1},x_{2}\right)  =0,~x_{2}\neq0
\end{array}
\end{equation}

\begin{figure}
\begin{center}
\includegraphics[scale=0.75]{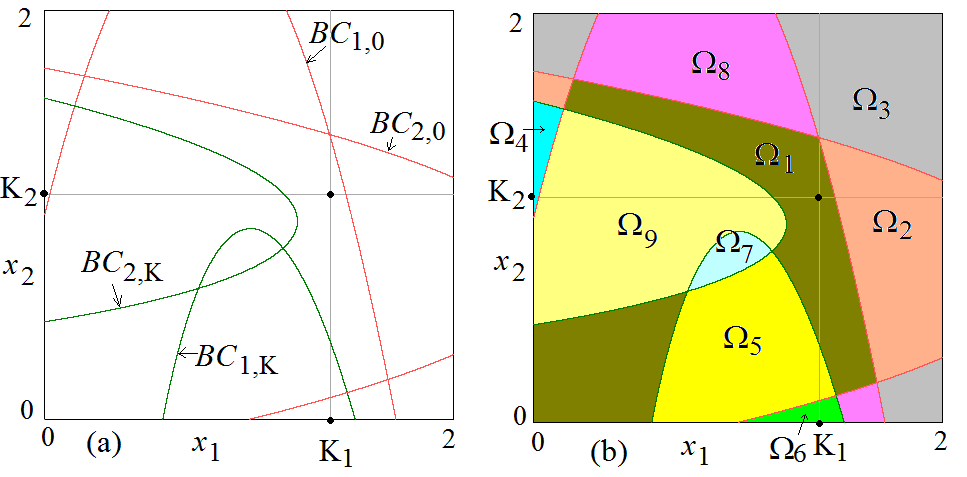}
\end{center}
\caption{Parameters: $\gamma_{1}=\gamma_{2}=1$, $\tau_{1}=\tau_{2}=4$,
$N_{1}=N_{2}=1.5,$ as in all the figures of the paper. $K_{1}=1.4$ and
$K_{2}=1.1$. The gray lines mark the border of the phase plane $D=\left[
0,K_{1}\right]  \times\left[  0,K_{2}\right]  $. A black dot marks the
intersection of these two lines. In (a) the border curves $BC_{1,0}$ and
$BC_{2,0}$ in red, the border curves $BC_{1,K}$ and $BC_{2,K}$ in green. In
(b) the regions $\Omega_{j}$ of the phase plane are evidenced by different
colors. To better clarify the shape of the constraints, we illustrate all the
figures in the phase plane $[0,2]\times\left[  0,2\right]  ,$ although the
region of interest of the model is $D\subset\left[  0,N_{1}\right]
\times\left[  0,N_{2}\right] $.}\label{Fig1}
\end{figure}

In Fig. \ref{Fig1}a these four curves are shown for parameter values $\gamma_{1}%
=\gamma_{2}=1$, $\tau_{1}=\tau_{2}=4$, $N_{1}=N_{2}=1.5$ and $K_{1}=1.4$ and
$K_{2}=1.1$. In the present work all the figures are shown with the values of
$\gamma_{i},$ $\tau_{i}$ and $N_{i}$ for $i=1,2$ as in Fig. \ref{Fig1}, while we let
vary the parameter values of $K_{1}$ and $K_{2},$ which are the constraints,
and are responsible for several border collision bifurcations.%

The positive quadrant of the phase plane is thus partitioned in nine regions,
in each of which a different definition (i.e. a different function) is to be
applied. Let us define the regions as follows:%
\begin{equation}%
\begin{array}
[c]{lll}%
\Omega_{1} & = & \left\{  \left(  x_{1},x_{2}\right)  |0\leq F_{1}\left(
x_{1},x_{2}\right)  \leq K_{1}\text{ and }0\leq F_{2}\left(  x_{1}%
,x_{2}\right)  \leq K_{2}\right\} \\
\Omega_{2} & = & \left\{  \left(  x_{1},x_{2}\right)  |F_{1}\left(
x_{1},x_{2}\right)  \leq0\text{ and }0\leq F_{2}\left(  x_{1},x_{2}\right)
\leq K_{2}\right\} \\
\Omega_{3} & = & \left\{  \left(  x_{1},x_{2}\right)  |F_{1}\left(
x_{1},x_{2}\right)  \leq0\text{ and }F_{2}\left(  x_{1},x_{2}\right)
\leq0\right\} \\
\Omega_{4} & = & \left\{  \left(  x_{1},x_{2}\right)  |F_{1}\left(
x_{1},x_{2}\right)  \leq0\text{ and }F_{2}\left(  x_{1},x_{2}\right)  \geq
K_{2}\right\} \\
\Omega_{5} & = & \left\{  \left(  x_{1},x_{2}\right)  |F_{1}\left(
x_{1},x_{2}\right)  \geq K_{1}\text{ and }0\leq F_{2}\left(  x_{1}%
,x_{2}\right)  \leq K_{2}\right\} \\
\Omega_{6} & = & \left\{  \left(  x_{1},x_{2}\right)  |F_{1}\left(
x_{1},x_{2}\right)  \geq K_{1}\text{ and }F_{2}\left(  x_{1},x_{2}\right)
\leq0\right\} \\
\Omega_{7} & = & \left\{  \left(  x_{1},x_{2}\right)  |F_{1}\left(
x_{1},x_{2}\right)  \geq K_{1}\text{ and }F_{2}\left(  x_{1},x_{2}\right)
\geq K_{2}\right\} \\
\Omega_{8} & = & \left\{  \left(  x_{1},x_{2}\right)  |0\leq F_{1}\left(
x_{1},x_{2}\right)  \leq K_{1}\text{ and }F_{2}\left(  x_{1},x_{2}\right)
\leq0\right\} \\
\Omega_{9} & = & \left\{  \left(  x_{1},x_{2}\right)  |0\leq F_{1}\left(
x_{1},x_{2}\right)  \leq K_{1}\text{ and }F_{2}\left(  x_{1},x_{2}\right)
\geq K_{2}\right\}
\end{array}
\end{equation}
so that the map in each region is given by:%
\begin{equation}%
\begin{array}
[c]{ll}%
\left(  x_{1},x_{2}\right)  \in\Omega_{1} & :(x_{1}^{\prime},x_{2}^{\prime
})=(F_{1}\left(  x_{1},x_{2}\right)  ,F_{2}\left(  x_{1},x_{2}\right)  )\\
\left(  x_{1},x_{2}\right)  \in\Omega_{2} & :(x_{1}^{\prime},x_{2}^{\prime
})=(0,F_{2}\left(  x_{1},x_{2}\right)  )\\
\left(  x_{1},x_{2}\right)  \in\Omega_{3} & :(x_{1}^{\prime},x_{2}^{\prime
})=(0,0)\\
\left(  x_{1},x_{2}\right)  \in\Omega_{4} & :(x_{1}^{\prime},x_{2}^{\prime
})=(0,K_{2})\\
\left(  x_{1},x_{2}\right)  \in\Omega_{5} & :(x_{1}^{\prime},x_{2}^{\prime
})=(K_{1},F_{2}\left(  x_{1},x_{2}\right)  )\\
\left(  x_{1},x_{2}\right)  \in\Omega_{6} & :(x_{1}^{\prime},x_{2}^{\prime
})=(K_{1},0)\\
\left(  x_{1},x_{2}\right)  \in\Omega_{7} & :(x_{1}^{\prime},x_{2}^{\prime
})=(K_{1},K_{2})\\
\left(  x_{1},x_{2}\right)  \in\Omega_{8} & :(x_{1}^{\prime},x_{2}^{\prime
})=(F_{1}\left(  x_{1},x_{2}\right)  ,0)\\
\left(  x_{1},x_{2}\right)  \in\Omega_{9} & :(x_{1}^{\prime},x_{2}^{\prime
})=(F_{1}\left(  x_{1},x_{2}\right)  ,K_{2})
\end{array}
\end{equation}
We notice that the points on the boundaries of the regions may belong to two
different regions: as the map is continuous, it does not matter whether a
point is considered belonging to one region or to the other, as the evaluated
value of the map is the same. From the definition it follows immediately that
the rectangle%
\begin{equation}
D=\left[  0,K_{1}\right]  \times\left[  0,K_{2}\right]
\end{equation}
is absorbing, as any point of the plane is mapped in $D$ in one iteration and
an orbit cannot escape from it, thus $D$ is our region of interest. In
general, depending on the values of the parameters, only a few of the regions
$\Omega_{j}$ for $j=1,...,9$\ may have a portion, or subregion, present in
$D$, say $\Omega_{j}\cap D\neq\varnothing,$ as shown for example in Fig. \ref{Fig1}b. In
any case, the behavior of the map in the other regions, not entering $D,$ may
be easily explained. To this purpose, let us introduce first a few remarks on
the fixed points that the system can have.

The fixed points of the system, satisfying $T(x_{1},x_{2})=(x_{1},x_{2})$, are
associated with the solutions of several equations. For sure we have some
fixed points on the axes, which correspond to disappearance (i.e. extinction)
of one population. From $F_{i}\left(  0,0\right)  =\left(  0,0\right)  $ for
$i=1,2$ we have that the origin $(0,0)$ is always a fixed point. Although, as
we shall see, it is locally unstable, all the points belonging to region
$\Omega_{3}$ are mapped into the origin in one iteration (and then they are fixed).

The axes are invariant, as considering a point $(x_{1},0)$ on the $x_{1}$ axis
we have that $T(x_{1},0)=(T_{1}(x_{1},0),T_{2}(x_{1},0))=(T_{1}(x_{1},0),0)$
still belongs to the axis, and
\begin{equation}
T_{1}(x_{1},0)=\left\{
\begin{array}
[c]{ccc}%
0 & \text{if} & F_{1}\left(  x_{1},0\right)  \leq0\\
F_{1}\left(  x_{1},0\right)  & \text{if} & 0\leq F_{1}\left(  x_{1},0\right)
\leq K_{1}\\
K_{1} & \text{if} & F_{1}\left(  x_{1},0\right)  \geq K_{1}%
\end{array}
\right.
\end{equation}
where
\begin{equation}
F_{1}\left(  x_{1},0\right)  =x_{1}\left[  1+\gamma_{1}x_{1}\tau_{1}\left(
1-\frac{x_{1}}{N_{1}}\right)  \right]  .
\end{equation}
Thus, region $\Omega_{6}$, whose points are all mapped in $(K_{1},0)$ in one
iteration, necessarily has non-empty intersection with the rectangle $D,$ and
$(K_{1},0)$ is a fixed point of the map. Moreover, considering the restriction
$t_{1}(x_{1})\equiv T_{1}\left(  x_{1},0\right)  $ we have that $t_{1}%
(x_{1})=x_{1}$ is satisfied for $x_{1}^{\ast}=0$ which is a fixed point
(representing the origin), and $x_{1}=N_{1}$\ which is virtual for
$K_{1}<N_{1}$ (constraint that we consider in the model). Thus the map
$x_{1}(t+1)=t_{1}(x_{1}(t))$ has a fixed point $x_{1}^{\ast}=K_{1}$ where the
piecewise smooth function $t_{1}(x_{1})$\ has a flat branch, which means that
the fixed point $(K_{1},0)$ always exists and is superstable (for the
restriction). While considering $\frac{d}{dx_{1}}t_{1}(x_{1})=1+2\gamma
_{1}\tau_{1}x_{1}-3\frac{\gamma_{1}\tau_{1}}{N_{1}}x_{1}^{2}$ we have
$\frac{d}{dx_{1}}t_{1}(0)=1$ and $\frac{d^{2}}{d^{2}x_{1}}t_{1}(x_{1}%
)=2\gamma_{1}\tau_{1}-6\frac{\gamma_{1}\tau_{1}}{N_{1}}x_{1}$ leads to
$\frac{d^{2}}{d^{2}x_{1}}t_{1}(0)=2\gamma_{1}\tau_{1}>0$ so that the fixed
point $x_{1}^{\ast}=0$ is repelling on its right side, that is, the origin
$(0,0)$ is repelling along the $x_{1}$ direction.

\begin{figure}
\begin{center}
\includegraphics[scale=0.75]{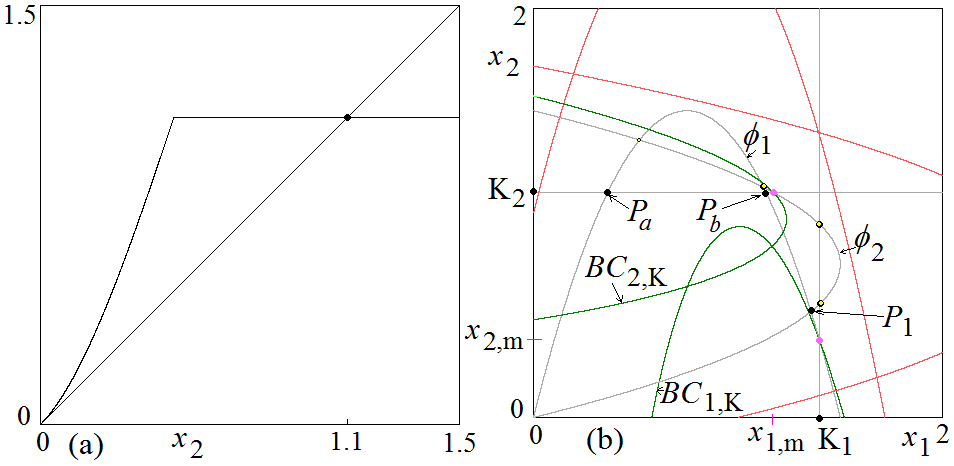}
\end{center}
\caption{Parameters as in Fig. \ref{Fig1}. In (a) function $t_{2}\left(  x_{2}\right)  $
for $x_{2}\in\left[  0,N_{2}\right]  $, the black dot is a superstable
equilibrium. In (b) the reaction curves $\phi_{1}$ and $\phi_{2},$ as well as
the lines $x_{1}=K_{1}$ and $x_{2}=K_{2},$ are in gray. Black dots are
feasible equilibria, black dots with yellow interior are virtual equilibria.
The dots in pink are the points $(x_{1,m},K_{2})$ and $(K_{1},x_{2,m}).$}\label{Fig2}
\end{figure}

Similarly for the second axis, we have that $T(0,x_{2})=(T_{1}(0,x_{2}%
),T_{2}(0,x_{2}))=(0,T_{2}(0,x_{2}))$ with%
\begin{equation}
T_{2}(0,x_{2})=\left\{
\begin{array}
[c]{ccc}%
0 & \text{if} & F_{2}\left(  0,x_{2}\right)  \leq0\\
F_{2}\left(  0,x_{2}\right)  & \text{if} & 0\leq F_{2}\left(  0,x_{2}\right)
\leq K_{2}\\
K_{2} & \text{if} & F_{2}\left(  0,x_{2}\right)  \geq K_{2}%
\end{array}
\right.
\end{equation}
where
\begin{equation}
F_{2}\left(  0,x_{2}\right)  =x_{2}\left[  1+\gamma_{2}x_{2}\tau_{2}\left(
1-\frac{x_{2}}{N_{2}}\right)  \right]  .
\end{equation}
So region $\Omega_{4}$ (whose points are all mapped in $(0,K_{2})$ in one
iteration) intersects the rectangle $D$ and $(0,K_{2})$ is a superstable fixed
point of the restriction, while the origin $(0,0)$ is repelling along the
$x_{2}$ direction. The proof is the same as the one given above for the
$x_{1}$ axis, changing the index $i=1$ into $i=2$. Regarding our example, the
one-dimensional map $x_{2}(t+1)=t_{2}(x_{2}(t))\equiv T_{2}(0,x_{2}(t))$ is
shown in Fig. \ref{Fig2}a.

Below we shall complete the comments regarding the fixed points $(K_{1},0)$
and $(0,K_{2})$\ on the axes for the two-dimensional map $T$.

\bigskip

Other fixed points $(x_{1}^{\ast},x_{2}^{\ast})$ may exist as solutions of the
equations%
\[
\left\{
\begin{array}
[c]{c}%
x_{1}R_{1}(x_{1})=x_{2}\\
x_{2}R_{2}(x_{2})=x_{1}%
\end{array}
\right.
\]
when belonging to region $\Omega_{1}$ (otherwise they are so-called virtual
fixed points). These fixed points can be seen in the phase plane as
intersection points of the two reaction curves%
\begin{equation}
\phi_{1}:\ x_{2}=x_{1}R_{1}(x_{1})\text{ and }\phi_{2}:x_{1}=x_{2}R_{2}%
(x_{2}),
\end{equation}
and the number of such points can be at most four.

Moreover, also fixed points $(K_{1},x_{2}^{\ast})$ may exist, associated with
the solutions of the equation $x_{2}R_{2}(x_{2})=K_{1},$ when belonging to
region $\Omega_{5}\cap D.$ Also these fixed points can be graphically seen in
the phase plane as intersection points of the two curves $x_{1}=K_{1} $
(vertical straight line) and $\phi_{2}$. Similarly, fixed points of type
$(x_{1}^{\ast},K_{2})$ associated with the solutions of the equation
$x_{1}R_{1}(x_{1})=K_{2}$ (intersection points of the horizontal straight line
$x_{2}=K_{2}$ and $\phi_{1}$) may exist, when belonging to region $\Omega
_{9}\cap D.$

The fixed points of the example shown in Fig. \ref{Fig1} are evidenced in Fig. \ref{Fig2}b where
the two curves$\ \phi_{1}$ and $\phi_{2}$ (having a unimodal shape) are drawn,
together with the straight lines $x_{1}=K_{1}$ and $x_{2}=K_{2}.$ Besides in
the origin, the curves$\ \phi_{1}$ and $\phi_{2}$ have three intersection
points, but two of them belong to region $\Omega_{9}$ and are outside $D$,
while the third one, say $P_{1},$ belongs to region $\Omega_{1} $ in $D$ and
thus it is a true fixed point of the map. On the vertical line $x_{1}=K_{1}$ a
fixed point is $(K_{1},0)$ on the axis and, as we shall see, it is
superstable. Then two more solutions of $x_{2}R_{2}(x_{2})=K_{1}$ exist, but
both points belong to region $\Omega_{1}$ and thus are virtual fixed points.
Differently, on the horizontal line\ $x_{2}=K_{2},$ besides the superstable
fixed point $(0,K_{2})$ on the vertical axis, there are two more fixed points
of the map, $P_{a}=(x_{1,a}^{\ast},K_{2})$ and $P_{b}=(x_{1,b}^{\ast},K_{2}),$
as both belong to region $\Omega_{9}\cap D.$ We shall return on these fixed
points below.

The definitions of the map in the several regions $\Omega_{j}$ lead to
different kinds of degeneracy. For example, when a portion of region
$\Omega_{7}$ exists in $D$, then all the points of that region are mapped into
a unique point: the corner $(K_{1},K_{2})$ of the absorbing rectangle $D,$
which means that in region $\Omega_{7}$ we have two degeneracies, that is, two
eigenvalues equal to zero in the Jacobian matrix at any point of $\Omega_{7}$.

Thus, one more fixed point may be given by the point $P=(K_{1},K_{2})$ when it
belongs to $\Omega_{7}\cap D$ (and in such a case this fixed point is
superstable: both eigenvalues are equal to zero). While when $(K_{1},K_{2})$
does not belong to $\Omega_{7}\cap D,$ then for the dynamics of the points in
region $\Omega_{7}$ it is enough to consider the trajectory of only one point:
$(K_{1},K_{2})$.

Other regions with double degeneracies are $\Omega_{3},$ $\Omega_{4}$ and
$\Omega_{6}$ as all of them are mapped into fixed points, $(0,0),$ $(0,K_{2})
$ and $(K_{1},0),$ respectively. These fixed points do not deserve for other
comments apart from their local stability/instability: as we have seen, the
origin is unstable while we shall see below (in Property \ref{Prop2}) that the two other
fixed points on the axes are superstable when $\Omega_{4}$ and $\Omega_{6}$
intersect $D$ in a set of positive measure, stable otherwise.\smallskip

There are other degeneracies which are immediate from the definition of the
map, due to the regions bounded by the border curves $BC_{i,K}$ (see Fig. \ref{Fig1}a).
Considering the portion of the phase plane which is bounded by the border
curve $BC_{1,K},$ we have that the whole region is mapped onto the line
$x_{1}=K_{1}.$ Similarly the whole region bounded by the border curve
$BC_{2,K}$ is mapped onto the line $x_{2}=K_{2}.$ Thus, in both regions we
have one degeneracy as the Jacobian matrix in all the points of these regions
has one eigenvalue equal to zero. As a whole region is mapped into a segment
of straight line, the dynamics can be associated with the points of those
particular segments. In particular, the stability/instability of the fixed
points belonging to these lines can be investigated considering the
restriction of the map to these lines, when they belong to the proper region
(that is, when they are real fixed points of $T$\ and not virtual). Let us
first notice the following\smallskip

\begin{property}\label{Prop1}
The three curves $x_{2}=K_{2}$, $BC_{2,K}$ and $\phi_{2}$ all intersect in the point $\left(x_{1,m},K_{2}\right)$, where $x_{1,m}=K_{2}\tau_{2}\left(1-\frac{K_{2}}{N_{2}}\right)$. The three curves $x_{1}=K_{1}$, $BC_{1,K}$ and $\phi_{1}$ all intersect in the point $\left(K_{1},x_{2,m}\right)$, where $x_{2,m}=K_{1}\tau_{1}\left(1-\frac{K_{1}}{N_{1}}\right)$.
\end{property}

\begin{proof}
In fact, $x_{2}=K_{2}$\ intersects $BC_{2,K}:x_{1}=\left[  1+\gamma_{2}%
x_{2}R_{2}\left(  x_{2}\right)  -\frac{K_{2}}{x_{2}}\right]  /\gamma_{2}$\ in
the point $x_{1,m}=K_{2}R_{2}(K_{2})=K_{2}\tau_{2}\left(  1-\frac{K_{2}}%
{N_{2}}\right)  ,$\ and also $x_{2}=K_{2}$\ intersects $\phi_{2}:x_{1}%
=x_{2}R_{2}(x_{2})=x_{2}\tau_{2}\left(  1-\frac{x_{2}}{N_{2}}\right)  $\ in
the same point, as it is immediately evident. Similarly for the other curves
(these points are evidenced in Fig. \ref{Fig2}b).
\end{proof}

So, let us consider $x_{2}=K_{2}$ and the segment of this line for $x_{1}%
\geq0$ and $x_{1}\leq x_{1,m},$ where $x_{1,m}$ is defined in Property \ref{Prop1}. Then
the restriction of the map to this segment is invariant, and on it the
dynamics are given (for $0\leq x_{1}\leq x_{1,m}$) by the one-dimensional map%
\begin{equation}
x_{1}(t+1)=f_{1}\left(  x_{1}(t)\right)  \ ,\ \ f_{1}\left(  x_{1}\right)
=\left\{
\begin{array}
[c]{ccc}%
0 & \text{if} & F_{1}\left(  x_{1},K_{2}\right)  \leq0\\
F_{1}\left(  x_{1},K_{2}\right)  & \text{if} & 0\leq F_{1}\left(  x_{1}%
,K_{2}\right)  \leq K_{1}\\
K_{1} & \text{if} & F_{1}\left(  x_{1},K_{2}\right)  \geq K_{1}%
\end{array}
\right. \label{restr}%
\end{equation}
where%
\begin{equation}
F_{1}\left(  x_{1},K_{2}\right)  =x_{1}\left[  1-\gamma_{1}K_{2}+\gamma
_{1}x_{1}\tau_{1}\left(  1-\frac{x_{1}}{N_{1}}\right)  \right] \label{F1K}%
\end{equation}
The point $x_{1}=0$ corresponds to the fixed point $(0,K_{2})$ of $T$, and
fixed points with positive values internal to the range $[0,K_{1}]$ are thus
associated with the solutions of a quadratic equation, leading to
\[
x_{1,b,a}^{\ast}=\frac{N_{1}}{2}\pm\sqrt{\left(  \frac{N_{1}}{2}\right)
^{2}-\frac{K_{2}N_{1}}{\tau_{1}}}
\]
moreover%
\begin{equation}
\frac{d}{dx_{1}}F_{1}\left(  x_{1},K_{2}\right)  =1-\gamma_{1}K_{2}%
+x_{1}\gamma_{1}\tau_{1}\left(  2-\frac{3x_{1}}{N_{1}}\right) \label{F1kder}%
\end{equation}
so that $\frac{d}{dx_{1}}F_{1}\left(  0,K_{2}\right)  =1-\gamma_{1}K_{2}<1$,
which implies that this fixed point is attracting also on the direction of the
line (as the derivative is either zero, when the constraint is active, or
positive and smaller that 1), and%
\begin{align*}
\frac{d}{dx_{1}}F_{1}\left(  x_{1,b,a}^{\ast},K_{2}\right)   & =1+2\gamma
_{1}K_{2}-x_{1,b,a}^{\ast}\gamma_{1}\tau_{1}\\
& =1+2\gamma_{1}K_{2}-\gamma_{1}\tau_{1}\left(  \frac{N_{1}}{2}\pm
\sqrt{\left(  \frac{N_{1}}{2}\right)  ^{2}-\frac{K_{2}N_{1}}{\tau_{1}}%
}\right)
\end{align*}
Summarizing, these two more are fixed points of the two-dimensional map only
if $x_{1,a}^{\ast}\leq x_{1,m}$ and $x_{1,b}^{\ast}\leq x_{1,m}$ (as it occurs
in the example shown in Fig. \ref{Fig2}b), and their stability depends on the value of
$\frac{d}{dx_{1}}f_{1}\left(  x_{1,b,a}\right)$. When $|\frac{d}{dx_{1}%
}f_{1}\left(  x_{1,b,a}\right)  |<1$ (resp. $>1$) the fixed points are
attracting (resp. repelling). In the example considered in Fig. \ref{Fig2}b both fixed
points $P_{a}=(x_{1,a}^{\ast},K_{2})$ and $P_{b}=(x_{1,b}^{\ast},K_{2})$ are repelling.

We can reason similarly for the restriction of the map on the straight line
$x_{1}=K_{1},$ for $0\leq x_{2}\leq x_{2,m},$ where $x_{2,m}$ is defined in
Property \ref{Prop1}, which is given by the one-dimensional map%
\begin{equation}
x_{2}(t+1)=f_{2}\left(  x_{2}(t)\right)  \ ,\ \ f_{2}\left(  x_{2}\right)
=\left\{
\begin{array}
[c]{ccc}%
0 & \text{if} & F_{2}\left(  K_{1},x_{2}\right)  \leq0\\
F_{2}\left(  K_{1},x_{2}\right)  & \text{if} & 0\leq F_{2}\left(  K_{1}%
,x_{2}\right)  \leq K_{2}\\
K_{2} & \text{if} & F_{2}\left(  K_{1},x_{2}\right)  \geq K_{2}%
\end{array}
\right. \label{T2k}%
\end{equation}
where
\begin{equation}
F_{2}\left(  K_{1},x_{2}\right)  =x_{2}\left[  1-\gamma_{2}K_{1}+\gamma
_{2}x_{2}\tau_{2}\left(  1-\frac{x_{2}}{N_{2}}\right)  \right] \label{F2k}%
\end{equation}

Thus, besides $x_{2}=0,$ which represents the fixed point ($K_{1},0$), the
fixed points are associated with the solutions of a quadratic equation,
leading to
\[
x_{2,b,a}^{\ast}=\frac{N_{2}}{2}\pm\sqrt{\left(  \frac{N_{2}}{2}\right)
^{2}-\frac{K_{1}N_{2}}{\tau_{2}}}
\]
moreover%
\begin{equation}
\frac{d}{dx_{2}}F_{2}\left(  K_{1},x_{2}\right)  =1-\gamma_{2}K_{1}%
+x_{2}\gamma_{2}\tau_{2}\left(  2-\frac{3x_{2}}{N_{2}}\right) \label{F2kder}%
\end{equation}
so that $\frac{d}{dx_{2}}F_{2}\left(  K_{1},0\right)  =1-\gamma_{2}K_{1}<1$,
which implies that this fixed point is attracting also on the direction of the
line, and%
\begin{align*}
\frac{d}{dx_{2}}F_{2}(K_{1},x_{2,b,a}^{\ast})  & =1+2\gamma_{2}K_{1}%
-x_{2,b,a}^{\ast}\gamma_{2}\tau_{2}\\
& =1+2\gamma_{2}K_{1}-\gamma_{2}\tau_{2}\left(  \frac{N_{2}}{2}\pm
\sqrt{\left(  \frac{N_{2}}{2}\right)  ^{2}-\frac{K_{1}N_{2}}{\tau_{2}}%
}\right)
\end{align*}
These solutions are fixed points of the two-dimensional map only if
$x_{2,a}^{\ast}\leq x_{2,m}$ and $x_{2,b}^{\ast}\leq x_{2,m}.$

With the parameter values used in the example shown in Fig. \ref{Fig2}b both the
inequalities given above are not satisfied and these points are called virtual
fixed points (i.e. they are not fixed points of the two-dimensional map).

We can also end the comments on the fixed points on the axes for the
two-dimensional map $T$. In fact, let us consider $\left(  0,K_{2}\right)  .$
We have already seen that along the axes $x_{1}=0$ there is a zero eigenvalue,
and now we can complete with the eigenvalue along the invariant segment on
$x_{2}=K_{2}.$ From the definition of the restriction in (\ref{restr}) and
(\ref{F1K}) we have that either this point (the origin of the restriction) is
superstable (which occurs when $\Omega_{4}$ intersects $D$ in a set of positive
measure), or stable, as we have\ $0<\frac{d}{dx_{1}}F_{1}\left(
0,K_{2}\right)  <1.$ Similarly we can reason for the other fixed point
$(K_{1},0).$

This leads to an important property of the model:\ the two single "segregation
states" always exist and attract some points of the phase plane. How many
points depends on the structure of the basins of attraction of these fixed
points, and on the existence or not of other attracting sets having states
with positive values (not converging to the axes). However, some results are
already known from the remarks written above: as all the points of the region
$\Omega_{2}$ are mapped into the $x_{2}$\ axis, which is trapping and on which
we know there is convergence to the fixed point $(0,K_{2})$, so we can
immediately conclude that all the points of region $\Omega_{2}$ belong to the
basin of attraction of $(0,K_{2}).$ Similarly, all the points of region
$\Omega_{8}$ belong to the basin of attraction of the fixed point $(K_{1},0).$
We shall see some examples below. We can so state the following\smallskip

\begin{property}\label{Prop2}
Two stable fixed points always exist in map $T$ given in
(\ref{T}): $\left(0,K_{2}\right)$ and $\left(K_{1},0\right)$. The points of
region $\Omega_{4}$ are mapped into $\left(0,K_{2}\right)$ and those
of region $\Omega_{2}$ converge to $\left(0,K_{2}\right)$. If $\Omega_{4}\cap D$ has positive measure, then $\left(0,K_{2}\right)$ is superstable for the two-dimensional map $T$. The points of region $\Omega_{6}$ are mapped into $\left(K_{1},0\right)$, and those of region $\Omega_{8}$ converge to
$\left(K_{1},0\right)$. If $\Omega_{6}\cap D$ has positive measure,
then $\left(K_{1},0\right)$ is superstable for the two-dimensional map $T$.
\end{property}

\begin{figure}
\begin{center}
\includegraphics[scale=0.75]{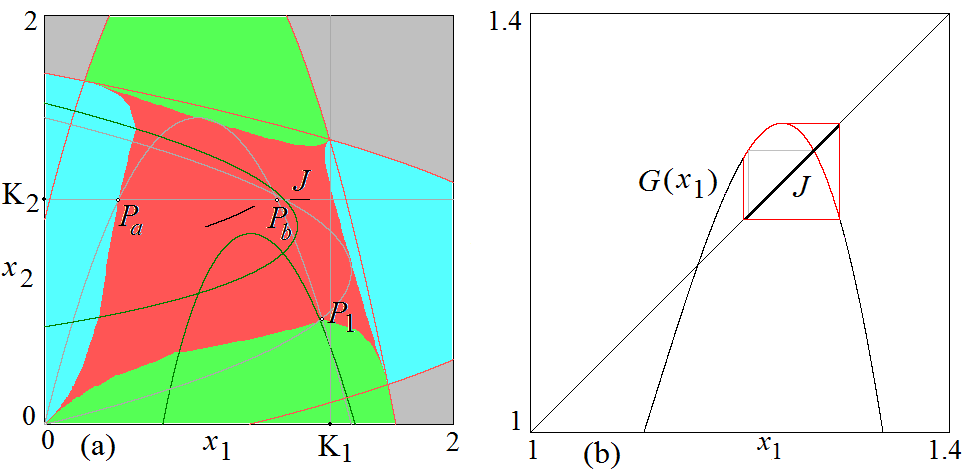}
\end{center}
\caption{Parameters as in Fig. \ref{Fig1}. In (a) basin of attraction of $\left(
0,K_{1}\right)  $ in green, basin of attraction of $\left(  0,K_{2}\right)  $
in azure, basin of attraction of $\left(  0,0\right)  $ in gray and basin of
attraction of the chaotic attractor in red. $P_{a}$, $P_{b}$ and $P_{1}$ are
unstable equilibria. In (b) first return map $G\left(  x_{1}\right)  $ on
$x_{2}=K_{2}$. $J$ represents an invariant segment, it is the portion of the
chaotic attractor of map $T$ that lies on $x_{2}=K_{2}$. The gray lines show
that the fixed point of the first return map is homoclinic.}\label{Fig3}
\end{figure}

In the example considered in Fig.s \ref{Fig1},\ref{Fig2} the two fixed points $(K_{1},0)$ and
$(0,K_{2})$\ on the axes are superstable for map $T$. Besides them, map $T$
has two more fixed points $P_{a}=(x_{1,a}^{\ast},K_{2})$ and $P_{b}%
=(x_{1,b}^{\ast},K_{2})$ in region $\Omega_{9}\cap D$ which are unstable, and
one more fixed point: $P_{1}\in\phi_{1}\cap\phi_{2}$ belonging to region
$\Omega_{1}\cap D.$ At $P_{1}$ the map has a smooth definition $(x_{1}%
(t+1),x_{2}(t+1))=(F_{1}(x_{1}(t),x_{2}(t)),F_{2}(x_{1}(t),x_{2}(t)))$, and
the stability of this fixed point depends on the eigenvalues of the Jacobian
matrix evaluated at $P_{1}.$ In our example also this fixed point $P_{1}$\ is
unstable. $P_{1}$ and $P_{a}$ belong to the frontiers separating the basins of
attraction. A third (chaotic) attractor exists, as shown in Fig. \ref{Fig3}a.

A trajectory on this attracting set consist of points which alternate from
region\ $\Omega_{9}$ to region $\Omega_{1}.$ This may be of great help as the
dynamics of $T$ can thus be investigated by use of a one dimensional map: the
first return map on a segment of the straight line $x_{2}=K_{2}$. In fact, the
points of the attracting set belonging to region $\Omega_{9}$ are mapped on
the line $x_{2}=K_{2}$ above the point $x_{1,m}=K_{2}\tau_{2}\left(
1-\frac{K_{2}}{N_{2}}\right)  $ (in region $\Omega_{1}).$\ Thus, a point
$(x_{1},K_{2})$ of the attractor is mapped in\ $T(x_{1},K_{2})=(F_{1}\left(
x_{1},K_{2}\right)  ,F_{2}\left(  x_{1},K_{2}\right)  )\in\Omega_{9}$ and then
a second iteration leads to $T^{2}(x_{1},K_{2})=(F_{1}(F_{1}\left(
x_{1},K_{2}\right)  ,F_{2}\left(  x_{1},K_{2}\right)  ),K_{2})=:(G\left(
x_{1}\right)  ,K_{2})\in\Omega_{1}.$ So it can be investigated by use of the
following one-dimensional first return map on $x_{2}=K_{2}$:%
\begin{equation}
x_{1}(t+1)=G\left(  x_{1}(t)\right)
\end{equation}%
\begin{equation}
G\left(  x_{1}\right)  =F_{1}\left(x_{1}\left(1-\gamma_{1}K_{2}+\gamma_{1}x_{1}\tau
_{1}\left(1-\frac{x_{1}}{N_{1}}\right)\right),K_{2}\left(1-\gamma_{2}x_{1}+\gamma_{2}K_{2}\tau
_{2}\left(1-\frac{K_{2}}{N_{2}}\right)\right)\right)
\end{equation}
in the range $x_{1,m}=K_{2}\tau_{2}\left(  1-\frac{K_{2}}{N_{2}}\right)
<x_{1}<K_{1}.$ This one-dimensional map, in our example, is shown in Fig. \ref{Fig3}b,
evidencing the invariant interval $J$ on which the dynamics seem to be
chaotic. Indeed, the fixed point in Fig. \ref{Fig3}b inside the invariant segment $J,$
which corresponds to an unstable 2-cycle of $T$, is homoclinic. This invariant
segment $J$ corresponds to the segment of the attractor on the straight line
$x_{2}=K_{2}$ in Fig. \ref{Fig3}a.\smallskip

As already remarked in the Introduction, the goal of this paper is to
investigate the role of the constraints, which are the values of $K_{1}$ and
$K_{2}.$ In doing so, here we investigate this only in the case in which the two states
(groups or populations) represented by $x_{1}$ and $x_{2}$ are in some way
symmetric, as characterized by parameters having the same values. Thus, in the
next section we shall consider the parameters $N\equiv N_{1}=N_{2}$,
$\tau\equiv\tau_{1}=\tau_{2}$ and $\gamma\equiv\gamma_{1}=\gamma_{2}$. Nevertheless, in piecewise smooth dynamical systems as the present one,
the other parameters may also be relevant. This aspect and in particular the
investigation of the role of the constraints in the generic case, with
different parameter values for the two populations, is left for further studies.

Here we are mainly interested in the role played by the two constraints
$K_{1}$ and $K_{2}$ which represent possible regulatory policy choices. Recall
that $K_{1}$ and $K_{2}$ represent the upper limit number of individuals of a
given group allowed to enter the system. We shall see a two-dimensional
bifurcation diagrams which immediately emphasizes the attracting cycles
existing as a function of the parameters $(K_{1},K_{2}).$ In the next section
we shall describe several regions in that parameter plane which lead to
interesting dynamic behaviors.

\section{Border Collision Bifurcations and global analysis of the dynamics}\label{BCBs}

Let us first consider the relevant dynamics occurring as a function of
$(K_{1},K_{2})$, let us call them {\it "the control parameters"}, when the other parameters are fixed (in our representative case at
the values considered in the figures of the previous section: $N=1.5$,
$\tau=4$ and $\gamma=1$). As in this paper we restrict our analysis to
populations with the same characteristics (in the parameters $\gamma_{i},$
$\tau_{i}$ and $N_{i}),$ the bifurcations occurring in the parameters
$(K_{1},K_{2})$ are obviously symmetric, which leads to the following
Property:

\begin{property}[Symmetric parameter plane]\label{Prop3}
Let $N\equiv N_{1}=N_{2}$, $\tau\equiv\tau_{1}=\tau_{2}$ and $\gamma\equiv\gamma_{1}=\gamma_{2}$. Let the control parameters have the
values $\left(K_{1},K_{2}\right)=\left(\xi,\eta\right)$ and let $\left\{
\left(a\left(t\right),b\left(t\right)\right), \; t>0\right\} $ be the trajectory associated
with the initial condition $\left(a\left(0\right),b\left(0\right)\right)$. Then $\left\{
\left(b\left(t\right),a\left(t\right)\right), \; t>0\right\}$ is the trajectory associated
with the initial condition $\left(b\left(0\right),a\left(0\right)\right)$ when the control parameters have the values $\left(K_{1},K_{2}\right)=\left(\eta,\xi\right)$.
\end{property}

\begin{figure}
\begin{center}
\includegraphics[scale=0.85]{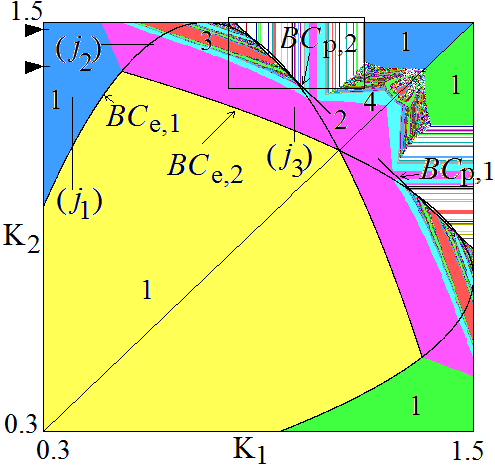}
\end{center}
\caption{Two dimensional bifurcation diagram in the $\left(  K_{1}%
,K_{2}\right)  $-parameter plane for the map $T$. Different colors are related
to attracting cycles of different periods $n\leq30$, the white region
corresponds either to chaotic attractors or to cycles of higher periodicity.}\label{Fig4}
\end{figure}

That is, via a change of variable $x_{2}:=x_{1}$ and $x_{1}:=x_{2}$ we have
the same dynamics when $K_{1}$ and $K_{2}$ are exchanged. This explains the
symmetric structure with respect to the main diagonal in the two-dimensional
bifurcation diagram shown in Fig. \ref{Fig4}.

As a particular case of Property \ref{Prop3} we have another property when $K_{1}=K_{2}
$ (on the diagonal of the two-dimensional bifurcation diagram of
Fig. \ref{Fig4}):

\begin{property}[Symmetric phase plane]
Let $N\equiv N_{1}=N_{2}$, $\tau\equiv\tau_{1}=\tau_{2}$, $\gamma\equiv
\gamma_{1}=\gamma_{2}$ and $K\equiv K_{1}=K_{2}$. Then:
\begin{enumerate}
\item[(4i)] Let $(a(t),b(t))$ for any integer $t>0$ be
the trajectory associated with the initial condition $\left(a\left(0\right),b\left(0\right)\right)$, then $\left(b\left(t\right),a\left(t\right)\right)$ for any integer $t>0$ is the
trajectory associated with the initial condition $\left(b\left(0\right),a\left(0\right)\right).$
\item[(4ii)] On the diagonal $\Delta$ of the phase plane map $T$ reduces to a one-dimensional system. From initial conditions $x_{1}\left(0\right)=x_{2}\left(0\right)$ it will be $x_{1}\left(t\right)=x_{2}\left(t\right)$ for any
integer $t>0$ and the iterates are given by the one-dimensional map
defined as $x\left(t+1\right)=T_{\Delta}\left(x\left(t\right)\right)$ with
\begin{equation}
T_{\Delta}(x)=\left\{
\begin{array}
[c]{ccc}%
0 & if & F_{\Delta}\left(  x\right)  \leq0\\
F_{\Delta}\left(  x\right)  & \text{if} & 0\leq F_{\Delta}\left(  x\right)
\leq K\\
K & if & F_{\Delta}\left(  x\right)  \geq K
\end{array}
\right. \label{Tdiag}%
\end{equation}
where $x\equiv x_{1}=x_{2}$ and $F_{\Delta}\left(
x\right)  \equiv F_{1}\left(  x,x\right)  =F_{2}\left(x,x\right)$ is given by%
\begin{equation}
F_{\Delta}\left(  x\right)  =x\left[  1-\gamma x+\gamma x\tau\left(
1-\frac{x}{N}\right)  \right]. \label{Fdiag}%
\end{equation}
\end{enumerate}
\end{property}

Clearly for the points of the phase plane outside the diagonal $x_{1}=x_{2}$
the Property (4i) stated above holds. Moreover, it is worth noting that Property (4i) implies that an invariant set of the two-dimensional
map is either symmetric with respect to the diagonal $\Delta$ $(x_{1}=x_{2}) $
of the phase plane, or the symmetric invariant set of it also exists.

As an example let us show the possible bifurcations occurring in the parameter
plane of the control parameters $(K_{1},K_{2})$ in the range $\left[
0,N_{1}\right]  \times\left[  0,N_{2}\right]  $ as reported in Fig. \ref{Fig4}. As the
model is symmetric (Property (4i)), we can just analyze the dynamics of the
model for $K_{2}\geq K_{1}$, i.e. taking into consideration only the region
above the diagonal in the two-dimensional bifurcation diagram of Fig. \ref{Fig4}, as the
dynamics and bifurcations for parameters on the symmetric side, i.e. for
$K_{2}\leq K_{1}$, are of the same kind (by Property \ref{Prop3}).

In Fig. \ref{Fig4} we highlights some BCB curves, which we shall explain below.

It is worth to note that as the parameters $K_{1}$ and $K_{2}$ influence the
borders of the regions at which the piecewise smooth map changes its
definition, all the bifurcations that we observe in Fig. \ref{Fig4} are expected to be
\textit{border collision bifurcations}. Indeed, even if this is not a
sufficient condition to state that all the curves are related to BCBs, the
high degeneracy of the map leads to this particular result.\smallskip

\subsection{Case $K_{1}=K_{2}$}\label{symcase}

Let us first describe the dynamics occurring in the phase plane when the
parameters belong to the diagonal $K_{1}=K_{2}$ of the two-dimensional
bifurcation diagram, and let $K\equiv K_{1}=K_{2}.$ As already shown above,
for points in the phase plane belonging to the diagonal where $x\equiv
x_{1}=x_{2}$ we can consider the one-dimensional piecewise smooth continuous
map $x(t+1)=T_{\Delta}(x(t))$ (given in (\ref{Tdiag}) and (\ref{Fdiag})).

The map $T_{\Delta}$ has fixed points satisfying the equation $F_{\Delta
}(x)=x$, leading to a fixed point in $x=0$ (representing the origin) and
$x^{\ast}=N(1-\frac{1}{\tau})$ which exists (positive) only for $\tau>1.$ It
is a real fixed point if $N(1-\frac{1}{\tau})\leq K,$ otherwise $x^{\ast}=K$
is a fixed point on the flat branch of the function. We can state the
following\smallskip

\begin{property}
Let $N\equiv N_{1}=N_{2}$, $\tau\equiv\tau_{1}=\tau_{2}>1$, $\gamma\equiv\gamma_{1}=\gamma_{2}$ and $K\equiv K_{1}=K_{2}$.
\begin{enumerate}
\item[(5i)] For $K<N\left(1-\frac{1}{\tau}\right)$ map $T_{\Delta}$ has a positive fixed point $x^{\ast}=K$ belonging to a flat branch, while for $K>N\left(1-\frac{1}{\tau}\right)$ map $T_{\Delta}$ has a positive fixed point $x^{\ast}=N(1-\frac{1}{\tau})$ belonging to a smooth branch. At $K=N\left(1-\frac{1}{\tau}\right)$ a border collision of the fixed point $x^{\ast}$ occurs.
If the bifurcation value satisfies $K<\overline{K}$ (resp. $K>\overline{K}$), where
\begin{equation}
\overline{K}=\frac{2\gamma(\tau-1)N+\sqrt{(2\gamma(\tau-1)N)^{2}+24N\gamma
\tau}}{6\gamma\tau}\label{eigen}%
\end{equation}
then increasing $K$ the result of the border collision is persistence of a stable fixed point (resp. a repelling fixed point and a superstable 2-cycle with periodic points $\left\{  K, \; T_{\Delta}\left(K\right)\right\}$).
\item[(5ii)] For $K>T_{\Delta}\left(x_{c}\right)$ where
\begin{equation}
x_{c}=\frac{(\tau-1)N}{3\tau}+\sqrt{\left(  \frac{(1-\tau)N}{3\tau}\right)  ^{2}+\frac{N}{3\gamma\tau}}\label{xc}%
\end{equation}
map $T_{\Delta}$ is smooth. At $K=T_{\Delta}(x_{c}$) there is a transition from piecewise-smooth to smooth.
\end{enumerate}
\end{property}

\begin{proof}
We notice that at $K=N(1-\frac{1}{\tau})$ for the two-dimensional map $T$ the
fixed point undergoes a codimension-two border collision as two borders are
crossed simultaneously $(\phi_{1}$ and $\phi_{2})$.

At the bifurcation value $K=N(1-\frac{1}{\tau})$ the fixed point $x^{\ast}$
merges with the border point (point in which the map changes its definition),
so it is a border collision. Increasing the value of $K$, the fixed point
$x^{\ast}$ moves from the flat branch to the smooth branch. The result of this
collision is completely predictable, as already remarked in the literature
(see for example \cite{SushkoGardini2010} and references therein). In fact, in the
one-dimensional case the skew-tent map can be used as a border collision
normal form, which means that in general, apart from codimension-two
bifurcation cases, the slopes of the two functions on the right and left side
of the border point at the BCB parameters values determine which kind of
dynamic behavior will appear after the BCB. In our case we have that the slope
on the left side of the border point is zero while on the right side it is
given by $F_{\Delta}^{\prime}\left(  K\right)  $ (also $F_{\Delta}^{\prime
}\left(  K\right)  =F_{\Delta}^{\prime}\left(  N(1-\frac{1}{\tau})\right)
=F_{\Delta}^{\prime}\left(  x^{\ast}\right)  ).$ Thus if $F_{\Delta}^{\prime
}\left(  K\right)  >-1$ (as the function is decreasing) we have persistence of
a stable fixed point, while if $F_{\Delta}^{\prime}\left(  K\right)  <-1$ the
fixed point on the smooth branch is unstable and a superstable 2-cycle exists
(i.e. with eigenvalue equal to zero). We have
\[
F_{\Delta}^{\prime}\left(  x\right)  =1+2\gamma(\tau-1)x-3\frac{\gamma\tau}%
{N}x^{2}
\]
so that $F_{\Delta}^{\prime}\left(  0\right)  =1$ and $F_{\Delta}%
^{\prime\prime}\left(  0\right)  =2\gamma(\tau-1)>0$ for $\tau>1$ leading to
$x^{\ast}=0$ repelling on its right side. Moreover,
\[
F_{\Delta}^{\prime}\left(  K\right)  =1+2\gamma(\tau-1)K-3\frac{\gamma\tau}%
{N}K^{2}
\]
and we have $F_{\Delta}^{\prime}\left(  K\right)  <-1$ for $K>\overline{K}$
where $\overline{K}$ is given in (\ref{eigen}), and for $K>N(1-\frac{1}{\tau
})$ a superstable 2-cycle appears, with periodic points $K$ and $F_{\Delta
}\left(  K\right)  $.

In our specific example considered in Fig. \ref{Fig4} the qualitative shape of the map
is shown in Fig. \ref{Fig5}a, it is $N(1-\frac{1}{\tau})=1.125,$ thus for $K<1.125$ the
map has a positive fixed point $x^{\ast}=K.$ The BCB of the fixed point occurs
at $K=1.125,$ and it is $\overline{K}\simeq1.07,$ so that at the bifurcation
value we have $K>\overline{K}$ and by Property (5i) a 2-cycle appears.

\begin{figure}
\begin{center}
\includegraphics[scale=0.55]{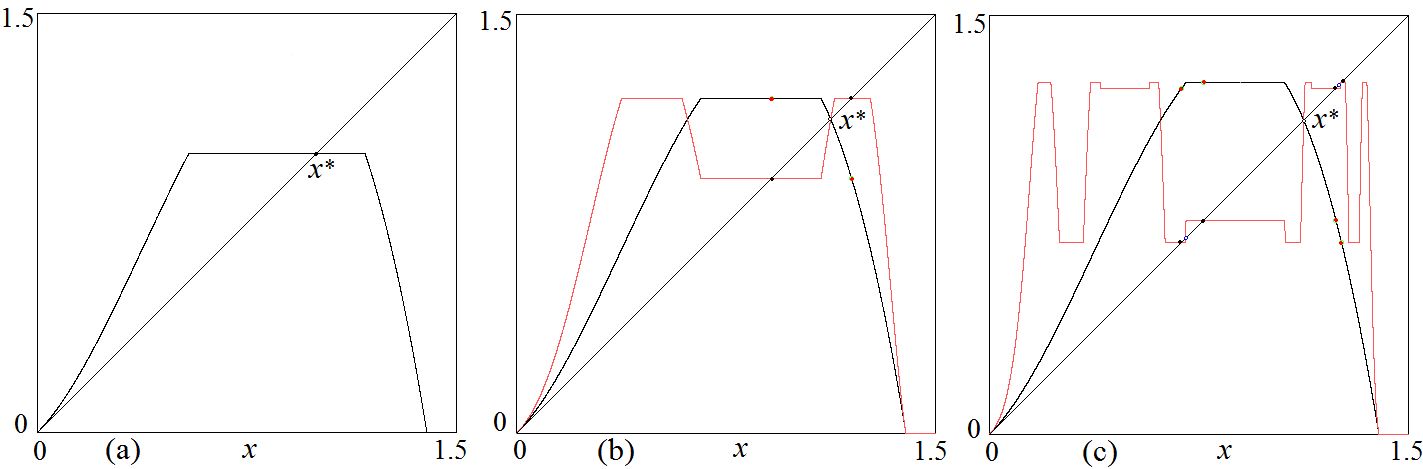}
\end{center}
\caption{ In (a) map $T_{\Delta}$\ at $K=1$, superstable fixed point $x^{\ast
}=K$. In (b) map $T_{\Delta}$ and its second iterate at $K=1.2$, superstable
2-cycle and unstable fixed point $x^{\ast}$. In (c) map $T_{\Delta}$ and its
fourth iterate at $K=1.26$, superstable 4-cycle and unstable fixed point
$x^{\ast}$.}\label{Fig5}
\end{figure}

When the fixed point $x^{\ast}$ exists, belonging to the decreasing branch (i.e.
after the border collision), from piecewise smooth the map may become smooth.
To detect this transition let us consider the critical point $x_{c}$ of
$T_{\Delta}$ (point in which the derivative of $F_{\Delta}$ in \eqref{Fdiag} vanishes), where $x_{c}$ is given in (\ref{xc}). Then for $K<T_{\Delta
}(x_{c})$ the map $T_{\Delta}$ has a horizontal flat branch (as it occurs in
our example in Fig. \ref{Fig5}), while for $K\geq T_{\Delta}(x_{c})$ the map is smooth
(as it occurs in our example for $K=1.4)$.
\end{proof}

Notice that the two border points of the map $T_{\Delta}(x),$ bounding the
flat branch, are given by the solutions of the equation $F_{\Delta}\left(
x\right)  =K,$ that is%
\[
x\left[  1-\gamma x+\gamma x\tau\left(  1-\frac{x}{N}\right)  \right]  =K
\]
As long as the fixed point $x^{\ast}=K$ exists in the flat branch, the two
border points are one smaller and one larger than $K$, while after its BCB
(with the largest border point) the two border points are both smaller than
$K$ (see Fig. \ref{Fig5}).

After the BCB of the fixed point we can consider the second iterate of the map
$T_{\Delta}^{2}(x)$ which,\ besides the unstable fixed point $x^{\ast
}=N(1-\frac{1}{\tau})$, has a pair of superstable fixed points (related to the
2-cycle) which also undergo a border collision. The BCB of the fixed point of
$T_{\Delta}^{2}$ can be studied in the same way as above for the fixed point
of $T_{\Delta}.$ In particular, a sequence of period doubling BCBs (also called flip BCBs)
occurs, leading to superstable cycles of period $2^{n}$.

In Fig. \ref{Fig5}a the fixed point $x^{\ast}=K$ is still on the flat branch, while in
Fig. \ref{Fig5}b, after its BCB, we have a 2-cycle, and in Fig. \ref{Fig5}c also the 2-cycle is
unstable and a superstable 4-cycle exists, with periodic points $K$\ and its
first three iterates.

As $K$ increases, all the cycles existing in the complete U-sequence (see
\cite{MetropolisSteinStein1973} and \cite{Hao1989}) appear also here, either by saddle-node BCB or by
flip BCB. For $T_{\Delta}(x)$ the cycles are either superstable or unstable.
The superstable cycles occur as long as in the map a flat branch persists,
that is, as remarked above in Property (5ii), as long as $K<T_{\Delta}%
(x_{c}),$ in which case the unstable cycles may belong to a chaotic repeller.
While for $K>T_{\Delta}(x_{c})$ an invariant chaotic set may exist for the
one-dimensional map $T_{\Delta}$ bounded by the critical point $T_{\Delta
}(x_{c})$ and its images. \smallskip

\begin{figure}
\begin{center}
\includegraphics[scale=0.75]{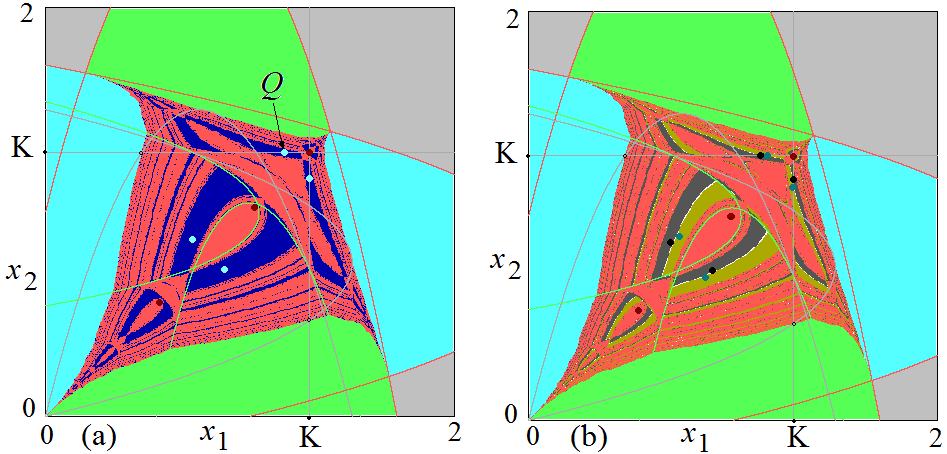}
\end{center}
\caption{Basin of attraction of $\left(  K,0\right)  $ in green. Basin of
attraction of $\left(  0,K\right)  $ in azure. Basin of attraction of $\left(
0,0\right)  $ in gray. Basin of attraction of the attractor lying on
$x_{1}=x_{2}$ in red. In (a) $K=1.2895,$ the blue region is the basin of
attraction of the 4-cycle. In (b) $K=1.29415,$ the dark-blue and dark-green
regions are the basins of attraction of the two 4-cycles born by pitchfork
bifurcation.}\label{Fig6}
\end{figure}

Going back to the two-dimensional map $T$ in the phase plane $(x_{1},x_{2})$,
for $K<T_{\Delta}(x_{c})$ the one-dimensional map $T_{\Delta}(x)$ is piecewise
smooth, and the attracting set for $T$ is some $n-$cycle on $\Delta$ having
one (and necessarily only one) periodic point belonging to region $\Omega_{7}$
and its image is the point $(K,K).$ \textit{It follows that such an }%
$n-$\textit{cycle is superstable also for the two-dimensional map }%
$T$\textit{.} However, it is not easy to predict the shape of its basin of
attraction, as this attractor coexists with the fixed points on the two axes,
and other attracting sets may exist in the phase plane outside $\Delta$. For
example, for $K=1.2$, when an attracting 2-cycle exists, its basin of
attraction is qualitatively similar to the one shown in Fig. \ref{Fig3}a for the chaotic
attractor. Differently it occurs for $K=1.2895$, when an attracting 3-cycle
exists on the diagonal $\Delta$, but it is not the unique attractor with
positive periodic points. In fact, it coexists with an attracting 4-cycle,
born in pair with an unstable 4-cycle via saddle-node BCB, and the stable set
of the unstable 4-cycle belongs to the frontier of the basins, shown in
Fig. \ref{Fig6}a.

In order to investigate the stability and bifurcations of the 4-cycle we
notice that, as already performed above, this can be done by use of a one
dimensional map: the first return map on the straight line$\footnote{We can
use, equivalently, the first return map on the straight line $x_{1}=K. $}$
$x_2=K$ for $x_1>x_{1,m}=K\tau\left(  1-\frac{K}{N}\right)$. So doing, it is
possible to consider $T^4(x_1,K)=(\psi(x_1),K)$ and the one-dimensional first
return map $x_1(t+1)=\psi(x_1(t))$ has a stable fixed point in the range
$[x_{1,m},K],$ corresponding to point $Q$ in Fig. \ref{Fig6}a, with a positive eigenvalue.
Increasing $K$\ this fixed point undergoes a pitchfork bifurcation, leading to
a pair of stable fixed points of $T^4$ which correspond to two stable 4-cycles
for $T$ (see Fig. \ref{Fig6}b). While the periodic points of the 4-cycles (one stable
and one unstable) in Fig. \ref{Fig6}a are symmetric with respect to $\Delta,$ those of
the pair of stable 4-cycles existing after the pitchfork bifurcation are not
symmetric themselves, but the two cycles have points which are pairwise
symmetric with respect to $\Delta$ (as stated in Property-(4i)).

\bigskip

\begin{remark}
Notice that even if we have called the described bifurcation
pitchfork, this term is proper only for the one-dimensional first return map
on the straight line $x_{2}=K.$\ In fact, let us reason as follows:
considering the attracting 4-cycle before the bifurcation (as shown in Fig. \ref{Fig6}a)
we can see that two periodic points are in region $\Omega_{1},$ one in region
$\Omega_{5}$ and one in region $\Omega_{9}$. Locally, in each point of the
4-cycle the map is smooth, and intuitively one can expect that the
stability/instability of the 4-cycle depends on the eigenvalues of the
Jacobian matrix of the map $T^{4}$ evaluated in any one of the four fixed
points belonging to the 4-cycle of $T$, and obviously one eigenvalue is
expected to be zero, due to the degeneracy of the map in regions $\Omega_{5}$
and $\Omega_{9}$. But this is not correct. The eigenvalue different from zero
so determined, is not associated with the bifurcations of the 4-cycle. This is
due to the degeneracy of the map: all the points of region $\Omega_{5}$ are
mapped onto the straight line $x_{1}=K$ independently on the eigenvalues
associated with the smooth map $T$ in points of this line belonging to region
$\Omega_{1}.$ That is: the bifurcation associated with cycles must be
determined by using the first return map, as we have done above, and not by
using the standard tools which are correct for smooth systems (also locally).
\end{remark}

\bigskip

\begin{figure}
\begin{center}
\includegraphics[scale=0.85]{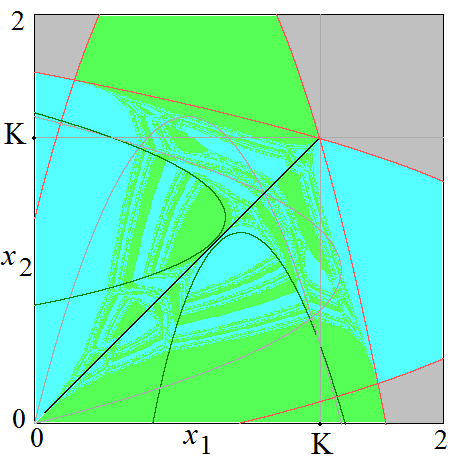}
\end{center}
\caption{ $K=1.4$ Basin of attraction of $\left(  K,0\right)  $ in green, basin
of attraction of $\left(  0,K\right)  $ in azure, separated by a fractal
frontier. The black dots along the line $x_{1}=x_{2}$ belong to a chaotic
saddle.}\label{Fig7}
\end{figure}

Differently from the case $K<T_{\Delta}(x_{c})$, when a superstable cycles
exists for $T$ on the diagonal of the phase plane, for $K>T_{\Delta}(x_{c})$
the one-dimensional map $T_{\Delta}(x)$ is smooth and an invariant set, which
may be chaotic, exists on $\Delta$ but this invariant set may be not
transversely attracting for the two-dimensional map $T$ in the phase plane
$(x_{1},x_{2}). $ In fact, this can also be observed in our example at $K=1.4$:
 a chaotic interval exists on the diagonal $\Delta$, which is a chaotic
repeller in the plane $\left(x_{1},x_{2}\right)$, the only attracting sets are the
fixed points $(K,0)$ and $(0,K)$ on the axes, and their basins are separated
by a fractal frontier, as shown in Fig. \ref{Fig7} (where the chaotic saddle is also
evidenced by a black segment on $\Delta$). This may lead to a significant
complexity in the socio-economic interpretation of the dynamics of the model. Indeed, given a generic value
$(x_{1}(0),x_{2}(0))$ as initial condition it is hard to predict whether the
states are ultimately converging to extinction of the first group or to
extinction of the second group.

The analysis conducted till now for $K_{1}=K_{2}$ reveals the importance of
the constraints for avoiding segregation. Indeed, from the dynamics of the
model we know that if the number of the members of the two populations that are allowed to
enter the system is sufficiently small, we always have a stable equilibrium of
non segregation. On the contrary, as the maximum number of agents of the two
groups that are allowed to enter the system increases, the equilibrium of non
segregation loses its stability and a sequence of cycles of different
periodicity appears. Further increasing this limit, we have that only
equilibria of segregation are stable. This positive effect of the entry
constraints on avoiding segregation can be explained observing that the
reaction of agents of one group toward the presence of agents of the opposed
group in the system is limited if the presence of the agents of both groups is
small in number. In other words, the entry constraints avoid the problem of
overshooting, which can be interpreted as impulsive and emotional behaviors.

\subsection{Case $K_{1}\neq K_{2}$}

\begin{figure}
\begin{center}
\includegraphics[scale=0.75]{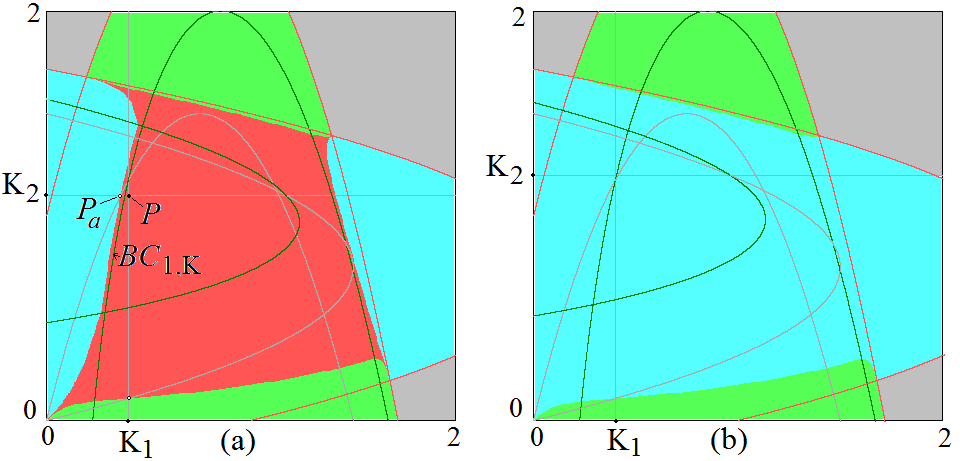}
\end{center}
\caption{Bifurcation through $(j_{1})$, i.e. $K_{1}=0.4$. In (a) a superstable
equilibrium, $P$, with basin in red, and a saddle, $P_{a}$, exist for
$K_{2}=1.1$. In (b) for $K_{2}=1.2$ the equilibria $P$ and $P_{a}$ do not
exist anymore as they disappeared by saddle-node BCB increasing $K_{2}$.}\label{Fig8}
\end{figure}

Let us first describe some of the BCB curves observable in Fig. \ref{Fig4}. The yellow
region in the center of the figure is associated with the existence of the
superstable fixed point $P=(K_{1},K_{2})\in\Omega_{7}.$ In our numerical
simulations (in the given example) it is the only attractor coexisting with
the fixed points on the axes, and its basin of attraction has a shape similar
to the one shown in Fig. \ref{Fig3}a (for the chaotic attractor). The boundaries of the
yellow region in the two-dimensional bifurcation diagram in Fig. \ref{Fig4} are clearly
curves of BCB, associated with a collision of $P$ with the borders $BC_{1,K}$
and $BC_{2,K}$ given in (\ref{BCK}). The condition for the border collision is
given by $P\in BC_{1,K}$ and $P\in BC_{2,K}$ leading to the BCB\ curves having
the following equations:%
\begin{align}
BC_{e,1}  & :K_{2}=K_{1}\tau_{1}\left(  1-\frac{K_{1}}{N_{1}}\right)  \text{
at which }P=(K_{1},K_{2})\in BC_{1,K}\\
BC_{e,2}  & :K_{1}=K_{2}\tau_{2}\left(  1-\frac{K_{2}}{N_{2}}\right)  \text{
at which }P=(K_{1},K_{2})\in BC_{2,K}\nonumber
\end{align}
which are drawn in Fig. \ref{Fig4}. Notice that the intersection point of these two BCB
curves, different from zero, is given by $K_{1}=K_{2}=N(1-\frac{1}{\tau})$
which corresponds to the BCB of the fixed point $K=N(1-\frac{1}{\tau})$
commented in Subsection \ref{symcase}. Let us consider the region with $K_{2}>K_{1}$, see Fig. \ref{Fig4}. For
parameters in the yellow region the fixed point $P$ is superstable. When a
parameter point crosses these curves the fixed point $P$ either disappear by
saddle-node BCB, when crossing $BC_{e,1},$ or enters (continuously) region
$\Omega_{5}$ when crossing $BC_{e,2}.$ In our example, for parameters crossing
$BC_{e,1}$ along the path $(j_{1})$ in Fig. \ref{Fig4}, the fixed point $P$ merges with
the unstable fixed point $P_{a}$ on the frontier of its basin of attraction
and disappears, leaving the two fixed points on the axes as the only
attractors. In Fig. \ref{Fig8}a it is shown the phase plane before the bifurcation, and
in Fig. \ref{Fig8}b after the bifurcation, when $P_{a}$ becomes virtual and
$(K_{1},K_{2})$ is no longer a fixed point. It can be seen that after the
bifurcation, the former basin of $P$ is included in the basin of $\left(
0,K_{2}\right)  .$

\begin{figure}
\begin{center}
\includegraphics[scale=0.55]{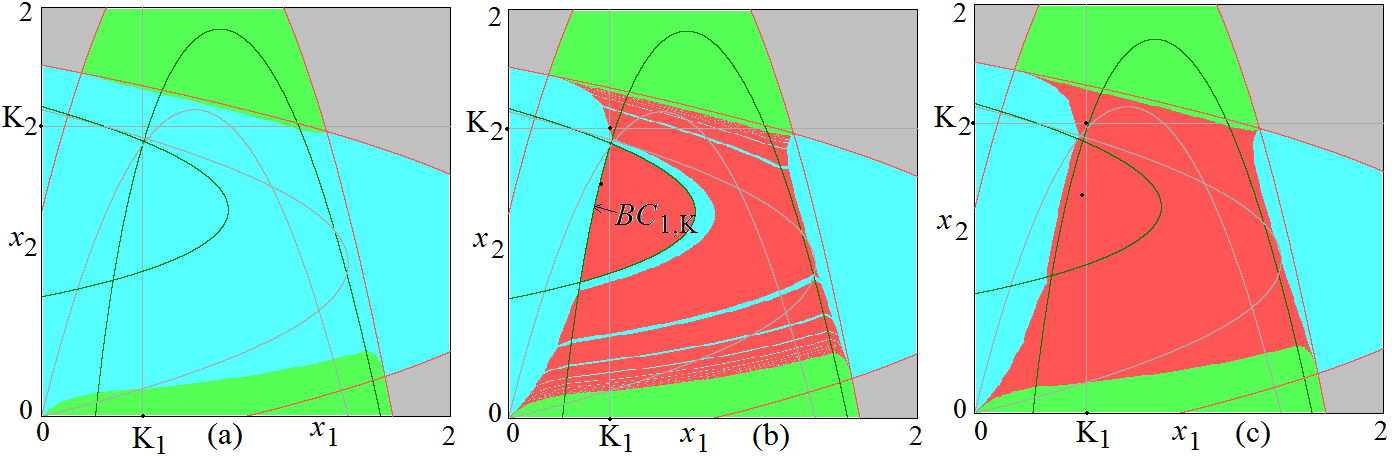}
\end{center}
\caption{Bifurcation through $(j_{2}),$ $K_{2}=1.42$. In (a) $K_{1}=0.4936.$ In
(b) $K_{1}=0.4937$. In (c) $K_{1}=0.55$. Basin of attraction of the 2-cycle in red.}\label{Fig9}
\end{figure}

A similar bifurcation involving a 2-cycle is shown changing the parameters
along the path $(j_{2})$ in Fig. \ref{Fig4}. For low values of $K_{1}$ only the two
fixed points on the axes are attracting (see Fig. \ref{Fig9}a). Increasing $K_{1},$ a
pair of 2-cycles appear by saddle-node BCB. Fig. \ref{Fig9}b shows the phase plane very
close to the bifurcation value, one of the pair of 2-cycles is attracting,
with one periodic point in region $\Omega_{7}$ and one in region $\Omega_{1}%
,$\ while the saddle 2-cycle has periodic points in regions $\Omega_{9}$ and
$\Omega_{1}$ (see Fig. \ref{Fig9}c).

The occurrence of this saddle-node BCB bifurcation of the 2-cycle can also be
determined analytically. In fact, considering the point $(K_{1},K_{2}),$ it
must be a fixed point for the second iterate of map $T$. Thus let
\begin{align}
F_{1}\left(  K_{1},K_{2}\right)   & =K_{1}\left[  1-\gamma_{1}K_{2}+\gamma
_{1}\tau_{1}K_{1}\left(  1-\frac{K_{1}}{N_{1}}\right)  \right] \\
F_{2}\left(  K_{1},K_{2}\right)   & =K_{2}\left[  1-\gamma_{2}K_{1}+\gamma
_{2}\tau_{2}K_{2}\left(  1-\frac{K_{2}}{N_{2}}\right)  \right] \nonumber
\end{align}
the BCB\ curve satisfies the equation%
\[
F_{1}\left(  F_{1}\left(  K_{1},K_{2}\right)  ,F_{2}\left(  K_{1}%
,K_{2}\right)  \right)  =K_{1}
\]
that is:%
\[
F_{1}\left(  K_{1},K_{2}\right)  \left[  1-\gamma_{1}F_{2}\left(  K_{1}%
,K_{2}\right)  +\gamma_{1}\tau_{1}F_{1}\left(  K_{1},K_{2}\right)  \left(
1-\frac{F_{1}\left(  K_{1},K_{2}\right)  }{N_{1}}\right)  \right]  =K_{1}
\]
Notice that in Fig. \ref{Fig4} we have plotted the complete curves $BC_{e,1}$ and
$BC_{e,2}$\ as also the other parts, not bounding the region of a superstable
fixed point, may be related to some border collision. Their effect may also be
only of {\it "persistence border collision"}, as it happens for example along the
path\ $(j_{2})$ in Fig. \ref{Fig4}: increasing $K_{1}$ the curve $BC_{e,1}$ is crossed,
and the stable 2-cycle persists stable, but with periodic points in different
regions (one point in $\Omega_{7}$ and one in region $\Omega_{9}).$%

\begin{figure}
\begin{center}
\includegraphics[scale=0.75]{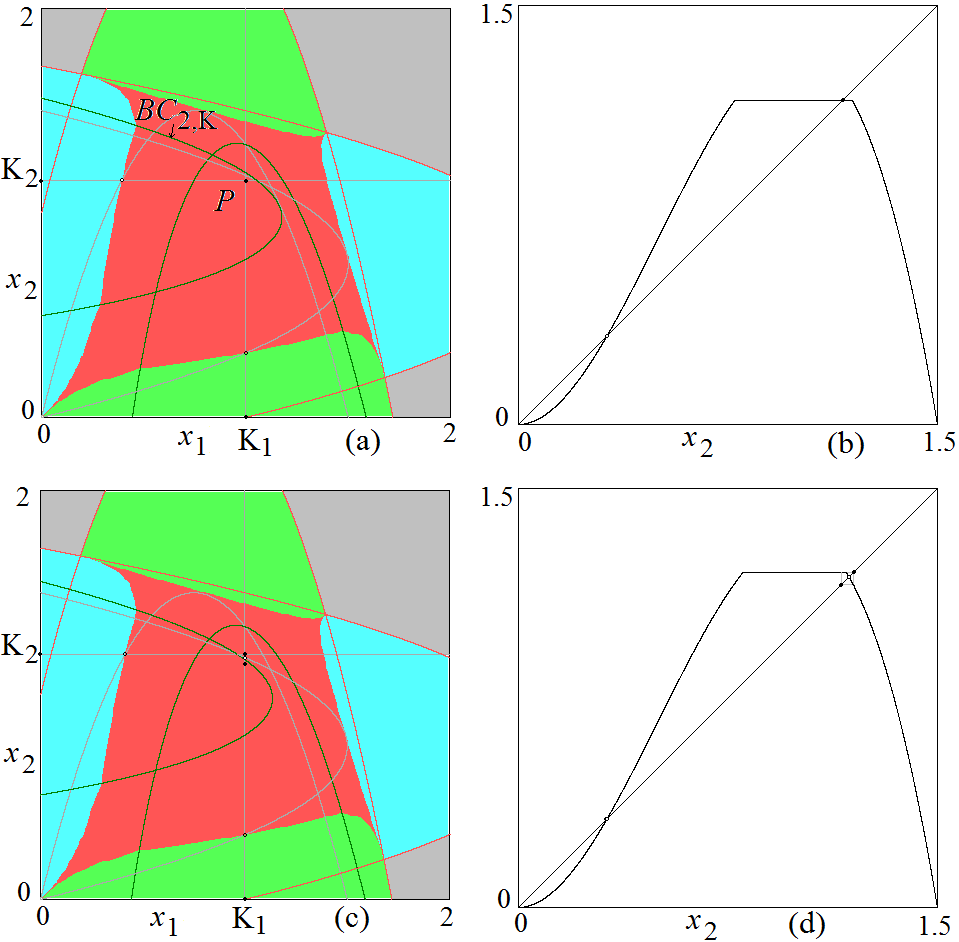}
\end{center}
\caption{Bifurcation through $(j_{3})$, i.e. $K_{1}=1.$ In (a) and (b)
$K_{2}=1.16$. In (a) the basin of attraction of $P=\left(  K_{1},K_{2}\right)
$ is in red. In (b) map $f_{2}\left(  x_{2}\right)  $ where black dot is the
superstable equilibrium $x_{2}^{\ast}=K_{2}$. In (c) and (d) $K_{2}=1.2$,
superstable 2-cycle appeared through a BCB of the fixed point $P$. In (c)
basin of attraction of this 2-cycle in red. In (d) map $f_{2}\left(
x_{2}\right)  $ where the black dots are the 2-cycle. In (a) and (c) the black
dots with while interior represent unstable equilibria.}\label{Fig10}
\end{figure}

In general, in order to predict the effect of the BCB of the fixed point, we
can use the first return map along the straight line $x_{1}=K_{1}$
(considering the part above the diagonal in Fig. \ref{Fig4}) and then make use of the
skew tent map as the border collision normal form, evaluating the slopes of
the functions at the border point, at the bifurcation values, as recalled in
the previous sections. For example, crossing the curve $BC_{e,2}$ along the
path $(j_{3})$ in Fig. \ref{Fig4} the fixed point $P$ crosses the curve $BC_{2,K}$ and
enters region $\Omega_{5}$. The fixed point becomes unstable and a stable
2-cycle appears, having one periodic point in region $\Omega_{7}$ and one in
region $\Omega_{5}$. In Fig. \ref{Fig10} it is shown the phase plane before the
bifurcation, and in Fig. \ref{Fig10}b the shape of the one-dimensional map restriction
of $T$ on the straight line $x_{1}=K_{1},$ for $0\leq x_{2}\leq x_{2,m},$
given in (\ref{T2k}) and (\ref{F2k}), showing the superstable fixed point on
the horizontal branch. The BCB of $P$ crossing the curve $BC_{2,K}$ in Fig. \ref{Fig10}a
corresponds to the BCB of the fixed point $x_{2}^{\ast}=K_{2}$ of the 1D map
(\ref{T2k}) in Fig. \ref{Fig10}b. The slopes at the bifurcation value are one zero and
one smaller than -1, thus the fixed point becomes unstable and a stable
2-cycle appears, as shown in Fig. \ref{Fig10}c,d. We can see that the structure of the
basins does not change.

\begin{figure}
\begin{center}
\includegraphics[scale=0.75]{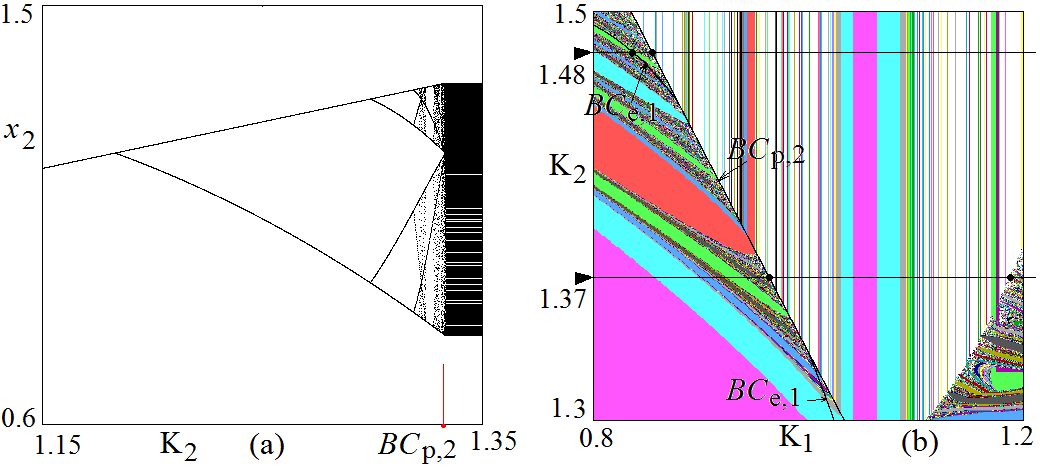}
\end{center}
\caption{In (a) 1D bifurcation diagram along $(j_{3})$, i.e. $K_{1}=1,$ and
$K_{2}\in\left[  1.15,1.35\right]  $. The red segment represents the moment in
which curve $BC_{p,2}$ is crossed. In (b) enlargement of the rectangle shown
in Fig. \ref{Fig4}.}\label{Fig11}
\end{figure}

Increasing $K_{2}$ along the path $(j_{3})$ in Fig. \ref{Fig4}, the one-dimensional
bifurcation diagram is reported in Fig. \ref{Fig11}a. It can be seen that after the
2-cycle, also attracting cycles of period 4 and $2^{n}$ for any $n$ exist.
This can be seen also in the enlargement of Fig. \ref{Fig4} reported in Fig. \ref{Fig11}b. This
region of the parameter plane corresponds to a region in which the BCBs lead
to the appearance of all stable cycles in accordance with the U-sequence, as
already remarked. In fact, the cycles there appearing all have one periodic
point in the region $\Omega_{7}$ and the periodic points either belong all to
the straight line $x_{1}=K_{1}$ (in which case the BCB can be studied via the
restriction of $T$\ on that line) or can be studied via the first return map
on that line. All these cycles are superstable for these one-dimensional maps
as well as for the two-dimensional map $T$, and undergo the border collisions.
The periodicity regions observable in Fig. \ref{Fig11}b are ordered according to the
U-sequence.

\begin{figure}
\begin{center}
\includegraphics[scale=0.75]{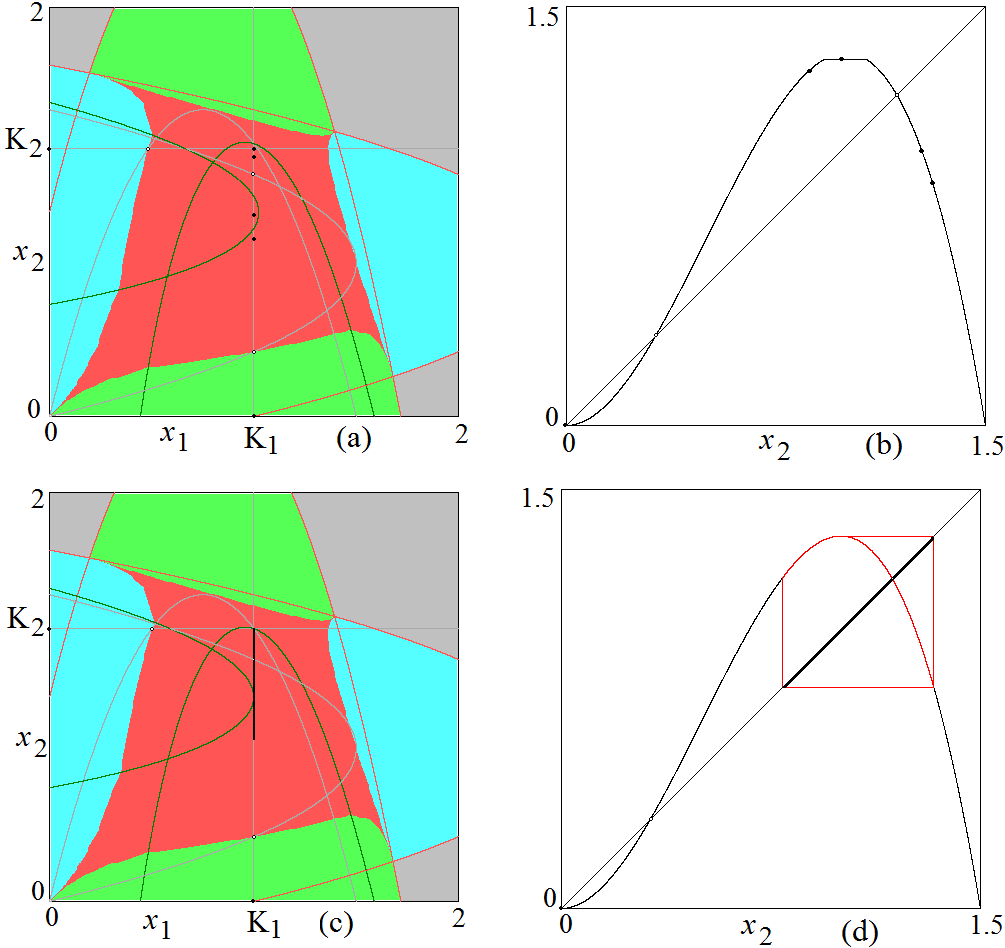}
\end{center}
\caption{$K_{1}=1$. In (a) and (b) $K_{2}=1.31$. In (a) basin of attraction of
the superstable 4-cycle (lying on $x_{1}=K_{1}$) in red. In (b) map
$f_{2}\left(  x_{2}\right)  $ describing the dynamics of the model on the
restriction $x_{1}=K_{1}$. In (c) and (d) $K_{2}=1.335$. In (c) basin of
attraction of the chaotic attractor (lying on $x_{1}=K_{1}$) in red. In (d)
map $f_{2}\left(  x_{2}\right)  $ which generates the chaotic attractor in
(c).}\label{Fig12}
\end{figure}

From the enlargement in Fig. \ref{Fig11}b it can be seen a change in the structure: the
periodicity regions of the superstable cycles (on the left side) end, and a
region with vertical strips appears. All the regions associated with
superstable cycles on the left, according to the U-sequence, also have
vertical strips on the right (still according with the U-sequence). This
transition, which is typical for one-dimensional piecewise smooth maps with a
horizontal branch, corresponds to the loss of the flat branch in the first
return map or in the one-dimensional restriction representing the dynamics of
the map $T$. In fact, as recalled above, the restriction of $T$\ on the line
$x_{1}=K_{1}$ has a horizontal branch as long as the cycles existing in the
region characterized by the U-sequence have one periodic point in region
$\Omega_{7}$. An example is shown in Fig. \ref{Fig10}, and in Fig. \ref{Fig12}a,b it is reported
the map at the value of $K_{1}$ for which there is a superstable 4-cycle.
Increasing $K_{1},$ a BCB occurs when the restriction of $T$\ to the line
$x_{1}=K_{1}$ becomes smooth, as shown in Fig. \ref{Fig12}c,d.

In order to obtain the bifurcation curves in the parameter space $(K_{1}%
,K_{2}),$ we proceed as follows. As recalled above, the restriction of map $T$
to the line $x_{1}=K_{1}$ is given in (\ref{T2k}) and (\ref{F2k}). The maximum
of the function $F_{2}\left(  K_{1},x_{2}\right)  ,$ $\max_{x_{2}}\left\{
F_{2}\left(  K_{1},x_{2}\right)  \right\}  $, is obtained considering its
proper critical point $x_{2,c}^{\ast},$ which satisfies $\frac{d}{dx_{2}}%
F_{2}\left(  K_{1},x_{2,c}^{\ast}\right)  =0,$ where the first derivative is
given in (\ref{F2kder}), and the value in the critical point (i.e.
$\max_{x_{2}}\left\{  F_{2}\left(  K_{1},x_{2}\right)  \right\}  =F_{2}\left(
K_{1},x_{2,c}^{\ast}\right)  )$. Via standard computations we get
\[
x_{2,c}^{\ast}=\frac{N_{2}}{3}+\left[  \left(  \frac{N_{2}}{3}\right)
^{2}+\frac{N_{2}}{3\gamma_{2}\tau_{2}}\left(  1-\gamma_{2}K_{1}\right)
\right]  ^{1/2}
\]
so that the maximum of the function $F_{2}\left(  K_{1},x_{2}\right)  $ is
given by
\[
F_{2}\left(  K_{1},x_{2,c}^{\ast}\right)  =\left(x_{2,c}^{\ast}\right)
^{2}\left(  2\frac{x_{2,c}^{\ast}}{N_{2}}-1\right)  \gamma_{2}\tau_{2}
\]
Then a BCB occurs when this maximum reaches the constraint on $x_{2},$ which
is the value $K_{2},$ and thus is determined by the condition $K_{2}%
=F_{2}\left(  K_{1},x_{2,c}^{\ast}\right)  $ which leads to the following BCB
curve in the parameter space:%
\begin{equation}
BC_{p,2}:K_{2}=\left(x_{2,c}^{\ast}\right)  ^{2}\left(  2\frac{x_{2,c}^{\ast}}{N_{2}}-1\right)  \gamma_{2}\tau_{2}\label{BP2}%
\end{equation}
A portion of this curve is shown in Fig. \ref{Fig4} and in the enlargement, in Fig. \ref{Fig11}b.

The other BCB due to the restriction on the straight line $x_{2}=K_{2}$ is
determined similarly, considering (\ref{restr}) and (\ref{F1K}). The maximum
of the function $F_{1}\left(  x_{1},K_{2}\right)  $ given in (\ref{F1K}) is
$\max_{x_{1}}\left\{  F_{1}\left(  x_{1},K_{2}\right)  \right\}  =F_{1}\left(
x_{1,c}^{\ast},K_{2}\right)  $, where $x_{1,c}^{\ast}$ is the proper critical
point, a solution of $\frac{d}{dx_{1}}F_{1}\left(  x_{1},K_{2}\right)  =0.$
From the first derivative given in (\ref{F1kder}) we get
\[
x_{1,c}^{\ast}=\frac{N_{1}}{3}+\left[  \left(  \frac{N_{1}}{3}\right)
^{2}+\frac{N_{1}}{3\gamma_{1}\tau_{1}}\left(  1-\gamma_{1}K_{2}\right)
\right]  ^{1/2}
\]
so that the maximum of the function is given by
\[
F_{1}\left(  x_{1,c}^{\ast},K_{2}\right)  =\left(x_{1,c}^{\ast}\right)
^{2}\left(  2\frac{x_{1,c}^{\ast}}{N_{1}}-1\right)  \gamma_{1}\tau_{1}
\]
A BCB occurs when this maximum reaches the value $K_{1},$ and thus is
determined by the condition $K_{1}=F_{1}\left(  x_{1,c}^{\ast},K_{2}\right)  $
which leads to the following BCB curve in the parameter space:%
\begin{equation}
BC_{p,1}:K_{1}=\left(x_{1,c}^{\ast}\right)  ^{2}\left(  2\frac{x_{1,c}^{\ast}%
}{N_{1}}-1\right)  \gamma_{1}\tau_{1}\label{BP1}%
\end{equation}

In Fig. \ref{Fig4} a portion of both bifurcation curves $BC_{p,1}$\ and $BC_{p,2}$\ are
shown, and better visible is $BC_{p,2}$ in the enlargement in Fig. \ref{Fig11}b. From
the two-dimensional bifurcation diagram we can see that the BCB occurring
crossing the curve $BC_{e,1}$ leads to persistence, while its portion in the
region with vertical strips is no longer a bifurcation, as the restriction to
the one-dimensional map is smooth and the point $(K_{1},K_{2})$ does not
belong to the attracting set.

Differently, the crossing of the curve $BC_{p,2}$ leading to a smooth
restriction, determines the transition from a piecewise-smooth (with a flat
branch) to a smooth map. It is worth to note that each periodicity region
associated with a superstable cycle, on the left side of the curve $BC_{p,2},
$ leads to a correspondent vertical strip associated with an attracting cycle
on its right side. On the left side of the curve $BC_{p,2}$ the periodicity
regions of superstable cycles have as limit sets curves related with
homoclinic bifurcations, which also leads to correspondent vertical lines
associated with chaotic dynamics on the right side (when the map is smooth).

In order to illustrate the dynamics of $T$ in this parameter region we
consider two more paths, at $K_{2}=1.48$ and at $K_{2}=1.37$ which also are
evidenced in Fig. \ref{Fig4} and Fig. \ref{Fig11}b, and describe some bifurcations occurring
as\ $K_{1}$ increases.

\begin{figure}
\begin{center}
\includegraphics[scale=0.73]{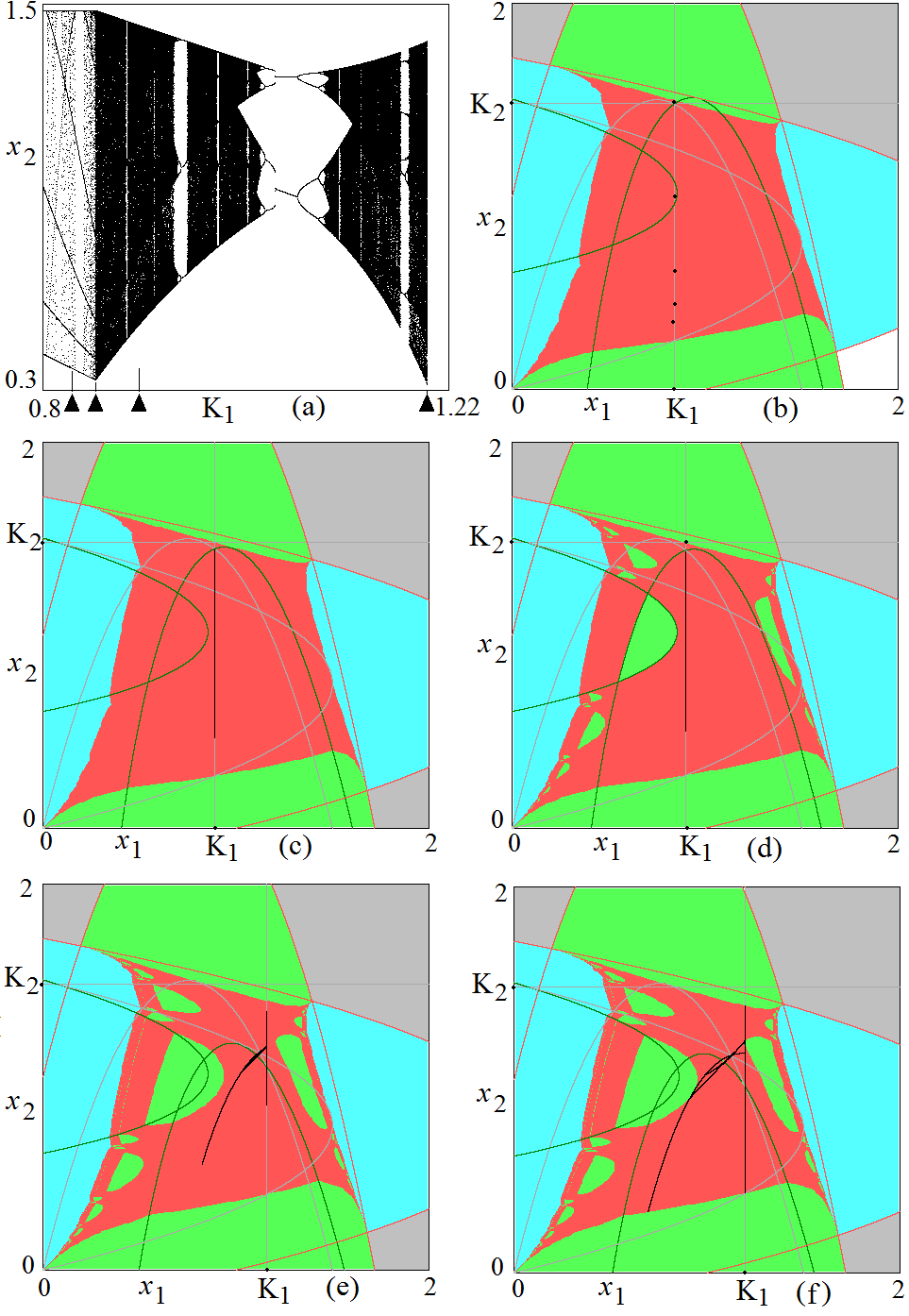}
\end{center}
\caption{$K_{2}=1.48$. In (a) 1D bifurcation diagram for $K_{1}\in\left[
0.8,1.22\right]  $. In (b) basins of attraction and attractors for $K_{1}=0.84
$. In (c) basins of attraction and attractors for $K_{1}=0.89$. In (d) basins
of attraction and attractors for $K_{1}=0.9$, here the corner point $\left(
K_{1},K_{2}\right)  $ is marked with a black dot for highlighting that it
enters the basin of attraction of $\left(  K_{1},0\right)  $, i.e. the green
region. In (e) basins of attractions and attractors for $K_{1}=1.16$. In (f)
basins of attraction and attractors for $K_{1}=1.198$.}\label{Fig13}
\end{figure}

Let us start considering $K_{2}=1.48$ fixed. From Fig. \ref{Fig11}b we can see that
increasing $K_{1}$ first the BCB crossing $BC_{e,1}$ occurs, and then the
crossing of $BC_{p,2}$. The one-dimensional bifurcation diagram as a function
of $K_{1}$ is shown in Fig. \ref{Fig13}a. In the region where the dynamics are
represented by the U-sequence as commented above, the effect of the crossing
of $BC_{e,1}$ (which occurs approximately at $K_{1}=0.83)$ corresponds to a
persistence of the attracting cycle: before the bifurcation the superstable
cycles have one periodic point in region $\Omega_{7}$ and all others in region
$\Omega_{5}$ while after the bifurcation the periodic point $(K_{1},K_{2})$
belongs to region $\Omega_{1},$ its preimage to region $\Omega_{7}$ and all
others in region $\Omega_{5}$. Then, increasing $K_{1}$ the crossing of
$BC_{p,2}$ occurs (approximately at $K_{1}=0.855),$ and this leads to a smooth
shape of the first return map on $x_{1}=K_{1}$ (as above in Fig. \ref{Fig12}d). The
attracting set on this line seems a large invariant chaotic interval, as shown
in Fig. \ref{Fig13}c. Notice that after this bifurcation, the corner point $(K_{1}%
,K_{2})$\ (belonging to region $\Omega_{1}$) does not belong to the attracting
set. This fact may lead to changes in the structure of the basins of
attraction of the attracting sets. As an example, in Fig. \ref{Fig13}c the corner point
is very close to the boundary separating the basin of the chaotic attractor
from the basin of the fixed point $(K_{1},0)$. In Fig. \ref{Fig13}d (at $K_{1}=0.9$) we
are at the contact:\ the corner point belongs to the boundary of the basin of
$(K_{1},0),$ as in fact the complete region $\Omega_{7}$ which is mapped into
$P$, now belongs to the basin of $(K_{1},0)$ together with all its preimages
of any rank. Increasing $K_{1}$ the attractor (a cycle or a chaotic attractor)
takes a more complex shape in the two-dimensional phase plane: the dynamics
can still be studied by using the first return map on the line $x_{1}=K_{1}$
but the number of points of a trajectory outside the line changes at each
iteration so that it is difficult to have it analytically, even in implicit
form (an example is shown in Fig. \ref{Fig13}e). In Fig. \ref{Fig11}b we also have evidenced the
point related to the "final bifurcation", as the positive attractor (here
chaotic) has a contact with the boundary of its basin of attraction (see
Fig. \ref{Fig13}f at $K_{1}=1.198$). After it is transformed into a chaotic repeller,
leaving only the two attracting fixed points on the axes, with a basin of
attraction having a fractal structure (similar to the one shown above in
Fig. \ref{Fig7}).\smallskip

\begin{figure}
\begin{center}
\includegraphics[scale=0.75]{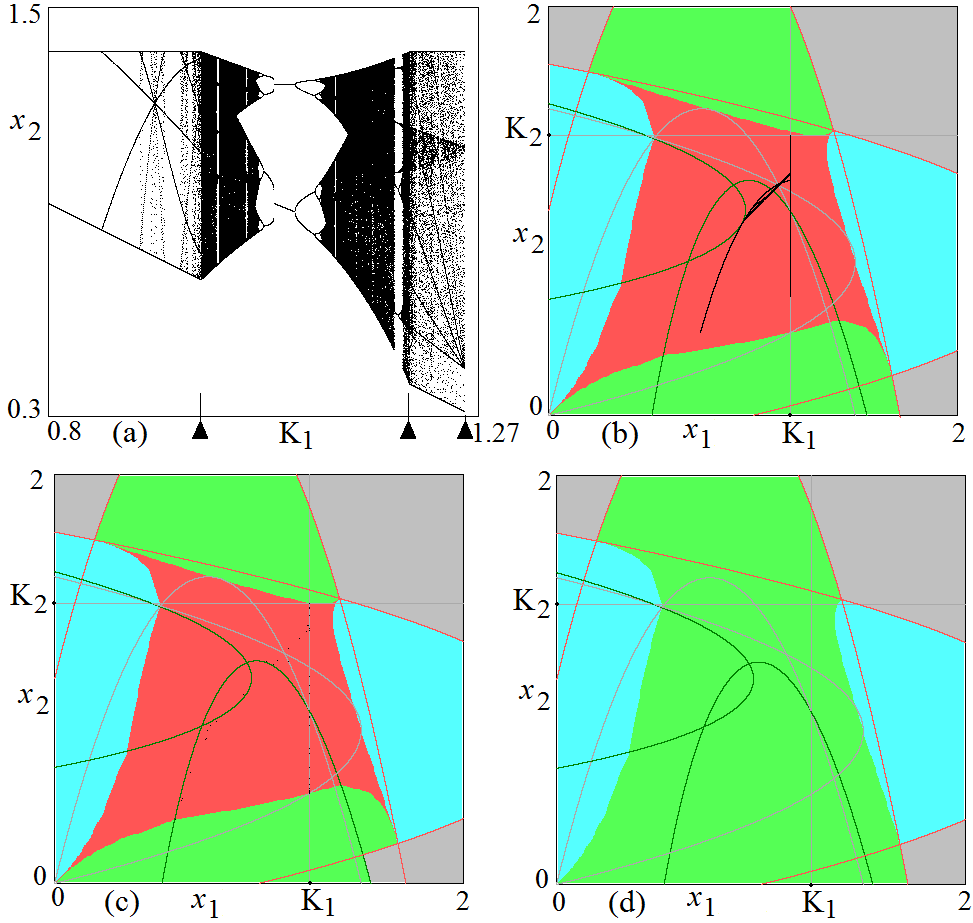}
\end{center}
\caption{$K_{2}=1.37$. In (a) 1D bifurcation diagram for $K_{1}\in\left[
0.8,1.27\right]  $. In (b) basins of attraction and attractors for
$K_{1}=1.185$. In (c) basins of attraction and attractors for $K_{1}=1.245$,
the cycle in the red region is superstable. In (d) basins of attraction and
attractors for $K_{1}=1.246$.}\label{Fig14}
\end{figure}

Let us now consider $K_{2}=1.37$ fixed. The one-dimensional bifurcation
diagram as a function of $K_{1}$ is shown in Fig. \ref{Fig14}a. From the region where
the dynamics are associated with superstable cycles and the U-sequence, the
BCB crossing $BC_{p,2}$ occurs approximately at $K_{1}=0.963,$ leading to a
chaotic attractor, which completely belongs to the line $x_{1}=K_{1}$ even if
after the bifurcation the point $(K_{1},K_{2})$ no longer belongs to the
attractor. The crossing of $BC_{e,1},$ which here occurs approximately at
$K_{1}=0.975,$ does not represent a bifurcation, it denotes only the
transition of the corner point $(K_{1},K_{2})$ from region $\Omega_{5}$ to
region $\Omega_{1}$. Increasing $K_{1}$ it can be noticed another region in
which the dynamics are again described by superstable cycles in the U-sequence
structure. This transition happens when the existing attractor has a contact,
i.e. a border collision, with the boundary of region $\Omega_{7}.$ In our
example this occurs approximately at $K_{1}=1.185$ as shown in Fig. \ref{Fig14}b. After
the contact the attractor is a superstable cycle with one periodic point in
region $\Omega_{7}$ and thus it is mapped into $P=(K_{1},K_{2})$ which is
again a periodic point, an example is shown in Fig. \ref{Fig14}c. The "final
bifurcation" of this attractor happens when the periodic point $P=(K_{1}%
,K_{2})$ has a contact with its basin boundary, which occurs approximately at
the value $K_{1}=1.245$ shown in Fig. \ref{Fig14}c:\ on the other side of the contact
point there is the basin of the fixed point $(K_{1},0$) so that after the
bifurcation the attractors are only the fixed points on the axes, and the
basin of\ $(K_{1},0$) increases, as shown in Fig. \ref{Fig14}d ($K_{1}=1.246$).%

For what concerns the implications of the entry constraints, $K_{1}$ and
$K_{2}$, in terms of segregation, the analysis conducted in this section
reveals that the effects of these entry constraints change if we make a strong
discrimination on the maximum number of agents allowed to enter the system
between the two groups. Indeed, if the difference between $K_{i}$ and $K_{j}$
is sufficiently large, with $K_{i}$ near to $N_{i}$ and $K_{j}$ small, then we
will have only stable equilibria of segregation. Moreover, starting with
$K_{i}$ relatively large and increasing $K_{j}$, a stable equilibrium of non
segregation cannot be reached, but rather an attractor in which the number of
agents of the two groups that enter and exit the system fluctuates over time
and when $K_{j}$ becomes sufficiently large again only an equilibrium of segregation
is possible. This reveals an important aspect of the issue of segregation,
i.e. to avoid overreaction of the two groups toward segregation we need to
limit in a similar way the number of possible entrances of both
types of agents in the system.

\section{Conclusions}\label{Conc}

In this work we have analyzed the effects of several constraints on the
dynamics of the adaptive model of segregation proposed in
\cite{BischiMerlone2011}. The constraints represent the maximum number of
agents of two different groups that are allowed to enter a system. We have
provided an accurate and deep investigation of the dynamics in the symmetric
case, i.e. when the two groups of agents that differ for a specific feature
are of the same size and have the same level of tolerance. The definition of
the two-dimensional piecewise smooth map lead to a map with different
definitions in several partitions. Besides the existence of two stable
segregation equilibria on the axes, we have shown that other attractors may
exist, regular or chaotic. The effect of the constraints, modifying the
regions, leads to border collision bifurcations of the positive attracting
sets. In the $\left(  K_{1},K_{2}\right) $-parameter plane of the
constraints, we have detected several BCB curves, explaining their effects on
the dynamic behaviors. The results are obtained by using several first return
maps on suitable intervals, and making use of the bifurcation theory
for one-dimensional piecewise smooth maps. A deep investigation of the effects
of the constraints when the symmetry is broken is desirable and can reveal
dynamics not observable in the symmetric setting. This line of research is
left for further studies.

\bigskip

\bigskip

\section*{ACKNOWLEDGMENTS}

This work has been performed under the activities of the Marie Curie
International Fellowship within the 7th European Community Framework
Programme, the project \textquotedblleft Multiple-discontinuity induced
bifurcations in theory and applications\textquotedblright. For the other two
authors also under the auspices of COST Action IS1104 "The EU in the new
complex geography of economic systems: models, tools and policy evaluation".

\bigskip

\bigskip

\bibliographystyle{plain}
%\bibliography{Thesis_ref}
\bibliography{C:/Users/UTENTE/Documents/PhD_Thesis_Davide_Radi/Thesis_ref}

\begin{thebibliography}{10}

\bibitem{AvrutinFutterSchanz2012}
V.~Avrutin, B.~Futter, and M.~Schanz.
\newblock The discontinuous top tent map and the nested period incrementing
  bifurcation structure.
\newblock {\em Chaos, Solitons \& Fractals}, 45:465--482, 2012.

\bibitem{BischiGardiniMerlone2009}
G.~I. Bischi, L.~Gardini, and U.~Merlone.
\newblock Impulsivity in binary choices and the emergence of periodicity.
\newblock {\em Discrete Dynamics in Nature and Society}, Volume 2009, 2009.

\bibitem{BischiMerlone2011}
G.~I. Bischi and U.~Merlone.
\newblock {\em Nonlinear economic dynamics}, chapter An Adaptive dynamic model
  of segregation, pages 191--205.
\newblock Nova Science Publisher, New York, 2011.

\bibitem{BrianzoniMichettiSushko2010}
S.~Brianzoni, E.~Michetti, and I.~Sushko.
\newblock Border collision bifurcations of superstable cycles in a
  one-dimensional piecewise smooth map.
\newblock {\em Mathematics and Computers in Simulation}, 81(1):52--61, 2010.

\bibitem{BischiChiarellaKopelSzidarovszky2009}
G.~I. Bischi~C. Chiarella, M.~Kopel, and F.~Szidarovszky.
\newblock {\em Nonlinear oligopolies: Stability and bifurcations}.
\newblock Heidelberg: Springer., 2009.

\bibitem{DalFornoGardiniMerlone2012}
A.~{Dal Forno}, L.~Gardini, and U.~Merlone.
\newblock Ternary choices in repeated games and border collision bifurcations.
\newblock {\em Chaos Solitons and Fractals}, in press doi:
  10.1016/j.chaos.2011.12.003, 2012.

\bibitem{Day1994}
R.~Day.
\newblock {\em Complex Economic Dynamics}.
\newblock MIT Press, Cambridge, 1994.

\bibitem{DayChen1993}
R.~Day and P.~Chen.
\newblock {\em Nonlinear Dynamics and Evolutionary Economics}, chapter
  Chaotically switching bear and bull markets: the derivation of stock price
  distributions from behavioral rules, pages 169--182.
\newblock Oxford University Press, Oxford, 1993.

\bibitem{Bernardo08}
M.~{di Bernardo}, C.~J. Budd, A.~R. Champneys, and P.~Kowalczyk.
\newblock {\em Piecewise-smooth Dynamical Systems: Theory and Applications}.
\newblock Springer-Verlag, Berlin, 2008.

\bibitem{EpsteinAxtell1996}
J.~M. Epstein and R.~L. Axtell.
\newblock {\em Growing Artificial Societies: Social Science from the Bottom
  up}.
\newblock Growing Artificial Societies: Social Science from the Bottom up.

\bibitem{GardiniMerloneTramontana2011}
L.~Gardini, U.~Merlone, and F.~Tramontana.
\newblock Inertia in binary choices: Continuity breaking and big-bang
  bifurcation points.
\newblock {\em Journal of Economic Behavior \& Organization}, 80(1):153--167,
  2011.

\bibitem{GardiniSushkoNaimzada2008}
L.~Gardini, I.~Sushko, and A.~Naimzada.
\newblock Growing through chaotic intervals.
\newblock {\em Journal of Economic Theory}, 143:541--557, 2008.

\bibitem{Hao1989}
B.~L. Hao.
\newblock {\em Elementary Symbolic Dynamics and Chaos in Dissipative Systems}.
\newblock World Scientific, Singapore, 1989.

\bibitem{HommesNusse1991}
C.~Hommes and H.~Nusse.
\newblock \textquotedblleft{P}eriod three to period two\textquotedblright
  bifurcation for piecewise linear models.
\newblock {\em Journal of Economics}, 54(2):157--169, 1991.

\bibitem{ItoTanakaNakada1979}
S.~Ito, S.~Tanaka, and H.~Nakada.
\newblock On unimodal transformations and chaos {II}.
\newblock {\em Tokyo Journal of Mathematics}, 2:241--259, 1979.

\bibitem{Matsuyama2010}
K.~Matsuyama.
\newblock The good, the bad, and the ugly: An inquiry into the causes and
  nature of credit cycles.
\newblock Center for Mathematical Studies in Economics and Management. Science
  Discussion Paper No.1391, Northwestern University., 2011.

\bibitem{MastrenkoMastrenkoChua1993}
Yu.~L. Ma\u{s}trenko, V.~L. Ma\u{s}trenko, and L.~O. Chua.
\newblock Cycles of chaotic intervals in a time-delayed {C}hua's circuit.
\newblock {\em International Journal of Bifurcation and Chaos in Applied
  Sciences and Engineering}, 3(6):1557--1572, 1993.

\bibitem{MetropolisSteinStein1973}
N.~Metropolis, M.~L. Stein, and P.~R. Stein.
\newblock On finite limit sets for transformations on the unit interval.
\newblock {\em J. Comb. Theory}, 15:25--44, 1973.

\bibitem{NusseYorke1992}
H.~E. Nusse and J.~A. Yorke.
\newblock Border-collision bifurcations including \textquotedblleft period two
  to period three\textquotedblright for piecewise smooth systems.
\newblock {\em Physica D}, 57(1--2):39--57, 1992.

\bibitem{NusseYorke1995}
H.~E. Nusse and J.~A. Yorke.
\newblock Border-collision bifurcations for piecewise smooth one-dimensional
  maps.
\newblock {\em International Journal of Bifurcation and Chaos in Applied
  Sciences and Engineering}, 5(1):189--207, 1995.

\bibitem{PuuSushko2002}
T.~Puu and I.~Sushko.
\newblock {\em Oligopoly Dynamics, Models and Tools}.
\newblock Springer Verlag, New York., 2002.

\bibitem{PuuSushko2006}
T.~Puu and I.~Sushko.
\newblock {\em Business Cycle Dynamics, Models and Tools}.
\newblock Springer Verlag, New York., 2006.

\bibitem{Schelling1969}
T.~C. Schelling.
\newblock Models of segregation.
\newblock {\em The American Economic Review}, 59(2):488--493, 1969.

\bibitem{SushkoAgliariGardini2006}
I.~Sushko, A.~Agliari, and L.~Gardini.
\newblock Bifurcation structure of parameter plane for a family of unimodal
  piecewise smooth maps: border-collision bifurcation curves.
\newblock {\em Chaos, Solitons \& Fractals}, 29(3):756--770, 2006.

\bibitem{SushkoGardini2010}
I.~Sushko and L.~Gardini.
\newblock Degenerate bifurcations and border collisions in piecewise smooth 1d
  and 2d maps.
\newblock {\em International Journal of Bifurcation and Chaos}, 20:2045--2070,
  2010.

\bibitem{SushkoGardiniMatsuyama2014}
I.~Sushko, L.~Gardini, and K.~Matsuyama.
\newblock \textquotedblleft superstable credit cycles and
  u-sequence\textquotedblright.
\newblock {\em Chaos Solitons \& Fractals}, 59:13--27, 2014.

\bibitem{TramontanaGardiniPuu2011}
F.~Tramontana, L.~Gardini, and T.~Puu.
\newblock Mathematical properties of a combined cournot-stackelberg model.
\newblock {\em Chaos, Solitons \& Fractals}, 44:58--70, 2011.

\bibitem{TramontanaGardiniWesterhoff2011}
F.~Tramontana, L.~Gardini, and F.~Westerhoff.
\newblock Heterogeneous {S}peculators and {A}sset {P}rice {D}ynamics: {F}urther
  {R}esults from a {O}n{e-D}imensional {D}iscontinuous {P}iecewis{e-L}inear
  {M}ap.
\newblock {\em Computational Economics}, 38(3):329--347, 2011.

\bibitem{TramontanaWesterhoffGardini2010}
F.~Tramontana, F.~Westerhoff, and L.~Gardini.
\newblock On the complicated price dynamics of a simple one-dimensional
  discontinuous financial market model with heterogeneous interacting traders.
\newblock {\em Journal of Economic Behavior \& Organization}, 74(3):187--205,
  2010.

\bibitem{Zhang2004}
J.~Zhang.
\newblock Residential segregation in an all-integrationist world.
\newblock {\em Journal of Economic Behavior and Organization}, 54:533--550,
  2004.

\bibitem{ZhusubaliyevMosekilde2003}
Z.~T. Zhusubaliyev and E.~Mosekilde.
\newblock {\em Bifurcations and Chaos in Piecewise-Smooth Dynamical Systems}.
\newblock World Scientific, River Edge, NJ, 2003.

\end{thebibliography}

\end{document}